\documentclass{article}
\usepackage[a4paper, margin=1in]{geometry}

\usepackage{amsmath}
\usepackage{nicefrac}
\usepackage{amssymb}
\usepackage{MnSymbol}
\usepackage{mathdots}
\usepackage{mathrsfs}
\usepackage{graphicx}
\usepackage{ctable}
\usepackage{cleveref}
\usepackage{siunitx}
\usepackage[utf8]{inputenc}

\DeclareMathOperator{\tridiag}{tridiag}
\DeclareMathOperator{\diag}{diag}
\DeclareMathOperator{\dx}{\Delta x}
\DeclareMathOperator{\dt}{\Delta t}
\DeclareMathOperator{\msym}{\mathscr{M}}
\DeclareMathOperator{\lsym}{\mathscr{L}}
\DeclareMathOperator{\gsym}{\mathscr{G}}
\DeclareMathOperator{\dsym}{\mathscr{D}}
\DeclareMathOperator{\esym}{\mathscr{E}}
\DeclareMathOperator{\ssym}{\mathscr{S}}

\newcommand{\Ogrande}{\mathcal{O}}

\newcommand{\C}{\mathbb{C}}

\newcommand{\p}{_{n}}
\newcommand{\glteq}{\sim_{\rm GLT}}

\newcommand{\dd}{\mathrm{d}}
\newcommand{\ii}{{\bf i}}

\providecommand{\keywords}[1]{\textbf{Keywords} #1}
\providecommand{\AMS}[1]{\textbf{MSC (2020)} #1}

\newtheorem{theorem}{Theorem}

\newtheorem{proposition}[theorem]{Proposition}
\newtheorem{definition}[theorem]{Definition}
\newtheorem{remark}[theorem]{Remark}

\newtheorem{corollary}[theorem]{Corollary}

\title{A matrix-theoretic spectral analysis of incompressible Navier-Stokes staggered DG approximation and related solvers}

\author{M. Mazza\thanks{
		Dipartimento di Scienze Umane e dell'Innovazione per il Territorio - 
		Universit\`a degli studi dell'Insubria. 
		Via Valleggio, 11 - 22100 Como.
		({mariarosa.mazza@uninsubria.it})
	} 
	\and M. Semplice\thanks{
		Dipartimento di Scienza e Alta Tecnologia - 
		Universit\`a degli studi dell'Insubria. 
		Via Valleggio, 11 - 22100 Como.
		({matteo.semplice@uninsubria.it})
	} 
	\and S. Serra Capizzano\thanks{
		Dipartimento di Scienze Umane e dell'Innovazione per il Territorio - 
		Universit\`a degli studi dell'Insubria. 
		Via Valleggio, 11 - 22100 Como.
		({s.serracapizzano@uninsubria.it})
	} 
	\and E. Travaglia \thanks{
		Dipartimento di Matematica -
		Universit\`a di Torino. 
		Via C. Alberto, 8 - 10123 Torino, Italy
		({elena.travaglia@unito.it})
	} 
}
\date{}

\begin{document}	
	\maketitle
	

\begin{abstract}
	The incompressible Navier-Stokes equations are solved in a channel, using a Discontinuous Galerkin method over staggered grids. The resulting linear systems are studied both in terms of the structure and in terms of the spectral features of the related coefficient matrices.  In fact, the resulting matrices are of block type, each block showing Toeplitz-like, band, and tensor structure at the same time. Using this rich matrix-theoretic information and the Toeplitz, Generalized Locally Toeplitz technology, a quite complete spectral analysis is presented, with the target of designing and analyzing fast iterative solvers for the associated large linear systems. 
	Quite promising numerical results are presented, commented, and critically discussed for elongated two- and three-dimensional geometries.
\end{abstract}

	\AMS{65F08, 65N30, 15B05, 15A18	}
	
	\keywords{Discontinous Galerkin,
		Incompressible Navier-Stokes,
		Schur complement,
		Toeplitz-like matrices,
		circulant preconditioning,
		spectral analysis}
 

\section{Introduction}

The efficient computation of incompressible fluid flows in complex geometries is a very important problem for physical and engineering applications.
In particular a delicate and time consuming task is the generation of the computational grid for a given geometry. 
Efficient algorithms avoid this step for example 
employing only a fixed background mesh
and discretizing
the equations for incompressible fluids with various strategies,
among which 
volume of fluid \cite{PR16:VOF},
ghost point \cite{Coco2020}, 
cut-cell \cite{CB10:LSSTAG:2d,NCB18:LSSTAG:3d,KK17}
and 
immersed boundary \cite{MittalIaccarino05:IBM} methods.
In all these methods, the description of the computational domain is often encoded in a level set function (see e.g. \cite{Sethian:LSbook,GibouFedkivOsher:reviewLevelSet}).
In particular, these techniques are very important in shape optimization problems since the mesh should be generated for all candidate geometries visited by the iterative optimization algorithm.

For industrial applications, a very important special case is the simulation of fluid flow in pipes of various cross-section. In this case, one can observe that the domain is much longer than wider and it is useful to leverage on one-dimensional or quasi-1D models, in which the pipe is described attaching a cross-section to each point of a 1D object. Notable examples in this direction are the Transversally Enriched Pipe Element Method of \cite{TEPEM:17} and the discretization methods at the base of the hierarchical model reduction techniques of \cite{HIModPipes:18}. Both of them compute a three-dimensional flow in a domain that is discretized only along the axial coordinate, i.e. the elements are sections of the whole pipe of length $\dx$. The finite element bases are obtained by Cartesian product of different discretizations in the longitudinal and in the transversal directions.

In this work we study a further simplification of the model, in which the transversal velocity components are neglected and only the longitudinal velocity is considered.
In particular we consider the incompressible Navier-Stokes equations
\begin{subequations}
\label{eq:INS}
\begin{gather}
\rho \left( \dfrac{\partial \textbf{u} }{\partial t}  + \nabla \cdot F_c \right) = - \nabla p + \nabla \cdot(\mu \ \nabla \textbf{u} ) \label{eq:INS:u}\\
\nabla \cdot \textbf{u}  = 0 \label{eq:INS:p},
\end{gather}
\end{subequations}
where $\textbf{x} = (x, y, z)$ is the vector of spatial coordinates and $t$ denotes the time, $p$ is the physical pressure and $\rho$ is the constant fluid density and $\mu$ is the viscosity which is a constant function if we consider a newtonian fluid. $F_c = \textbf{u}  \otimes \textbf{u} $ is the flux tensor of the nonlinear convective terms, $\textbf{u}  = (u,v,w)$ is the velocity vector where $u$ is the component parallel to the pipe axis, while $v$ and $w$ are the transversal ones.

We consider as domain a pipe with a variable cross-section and
since it has a length much greater than the section, we neglect the transverse velocities, i.e. we assume $ v = w = 0 $ (and consequently also $\partial_y p=\partial_z p=0$), but we consider the dependence on the three spatial variables of the longitudinal component, i.e. $ u = u (x, y, z) $.
The discretization is then performed with Discontinous Galerkin methods on a staggered grid arrangement, i.e. velocity elements are dual to the main grid of the pressure elements, similarly to \cite{Tavelli:14,Tavelli:15}, leading to a saddle point problem for the longitudinal velocity and the pressure variables. 

Having in mind the efficient solution of such linear system, in this paper we focus on the spectral study of the coefficient matrix as well as of its blocks and Schur complement. More specifically, we first recognize that all the matrix coefficient blocks show a block Generalized Locally Toeplitz (GLT) structure and that, as such, can be equipped with a symbol. Second, we leverage on the symbols of the blocks to retrieve the symbol of the Schur complement and the symbol of the coefficient matrix itself. We stress that in order to accomplish these goals, we introduce some new spectral tools that ease the symbol computation when rectangular matrices are involved. In this setting we can deliver a block circulant preconditioner for the Schur complement that provides a constant number of iterations as the matrix-size increases and that, once nested into a Krylov-type solver for the original coefficient matrix, brings to lower CPU timings when compared with other state-of-the-art strategies.

The paper is organized as follows. In \S\ref{sec:discr} we describe in details the discretization of the quasi-1D incompressible Navier-Stokes model; in \S\ref{sec:prelim} we both recall the Toeplitz and GLT technology and we introduce some new spectral tools that will be used in \S\ref{sec:spectr} to perform the spectral analysis of the matrix of the saddle point problem. This leads to the proposal of an efficient optimal preconditioner for our system, which is tested in the numerical section \S\ref{sec:numer}.

\section{Discretization}
\label{sec:discr}

\begin{figure}
	\begin{center}
		\includegraphics{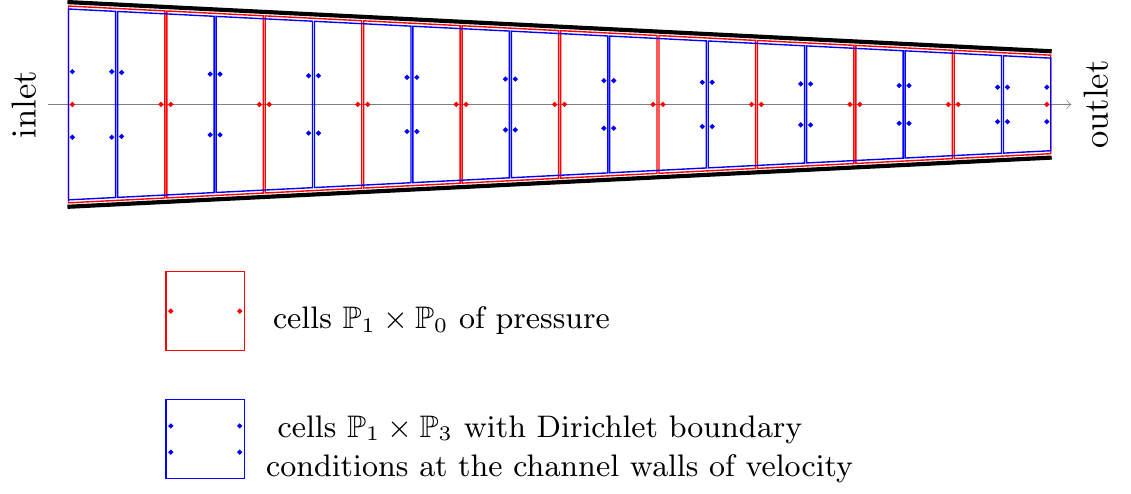}
	\end{center}
	\caption{Illustration of the staggered grid arrangement in a nozzle for $n_{x}=1$ and  $n_{y}=3$}
	\label{fig:discretiz}
\end{figure}

We consider the incompressible Navier-Stokes equations \eqref{eq:INS} in an elongated pipe-like domain, with a variable cross-section. An example is depicted in Fig.~\ref{fig:discretiz}.
We impose a no-slip condition at the solid boundaries;
at the outlet boundary we fix a null pressure, while at the inlet we impose Dirichlet data with a given velocity profile.

The channel is discretized only along its longitudinal dimension, so each cell is a section of the entire pipe of length $\dx$, (see Fig. \ref{fig:discretiz}). We denote the cells in this grid by $\Omega_{1}, \dots ,\Omega_{n}$. The discrete pressure is defined on this grid, while for the velocity we use a dual grid, whose first and last element have length equal to one half of the other cells. This type of staggered grid has been employed for example in \cite{Tavelli:14, Tavelli:15}. We denote the cells of the dual grid by $\Omega_{1}^{*}, \ldots ,\Omega_{n+1}^{*}$ and point out that each $\Omega_{j}$ has a nontrivial intersection only with $\Omega_{j}^{*}$ and $\Omega_{j+1}^{*}$ for $j= 1 \ldots n$.

For ease of presentation, we concentrate mainly on the two-dimensional case
and denote the width of the channel at the position $x$ by $d(x)$.
The longitudinal velocity $ u = u (x, y) $, in each cell of the dual grid,  is approximated by a $\mathbb{P}_{n_x}\otimes\mathbb{P}_{n_y}$ polynomial defined as the tensor product of the one dimensional polynomial of degree $n_x$ in the longitudinal direction and $n_y$ in the transverse one.
In order to do this, we construct a polynomial basis on the standard reference elements, $\Omega_{ref} = [0,1]^2$, using the Lagrange interpolation polynomials with equispaced nodes.
Taking into account the no-slip boundary condition applied at the channel walls, there are $n_u:=(n_x+1)\times(n_y-1)$ effective degrees of freedom for $u$ in each cell (blue dots in Fig.~\ref{fig:discretiz}). We stress that in order to satisfy the no-slip boundary conditions one should take $n_{y}\geq 2$.
\\In the same way the pressure is approximated in each cell of the primal grid by a $\mathbb{P}_{n_x}\otimes\mathbb{P}_{0}$ polynomial, i.e. the pressure is constant in the transversal direction. For this reason, there are only $n_p:=(n_x+1)$ degrees of freedom for $p$ in each cell (red dots in Fig.~\ref{fig:discretiz}). In general we are interested in a low degree $n_x$ but high degrees $n_y$, which are needed to compensate for the lack of mesh discretization in the transversal direction, and of course a mild but generic dependence of $d$ upon $x$.

To obtain a DG discretization on the staggered cell arrangements, 
we first integrate the momentum equation \eqref{eq:INS:u} multiplied by a generic shape function $\psi$ for the velocity over a cell of the dual grid, $\Omega_{i}^{*}$, for $ i = 1 \ldots n+1 $,
\begin{subequations}
\label{eq:weakINS}
\begin{equation}\label{eq:weakINS:u}
	\int_{\Omega_{i}^{*}}  \psi\ \rho \left( \dfrac{\partial u}{\partial t}  + \nabla \cdot F_c \right)d \textbf{x} 
	=
	- \int_{\Omega_{i}^{*}}  \psi\ \nabla p\  d \textbf{x} 
	+ \int_{\Omega_{i}^{*}}  \psi\ \nabla \cdot(\mu \ \nabla u) d \textbf{x}
.
\end{equation}
We then integrate the continuity equation \eqref{eq:INS:p},  multiplied by a generic shape function $\theta$ for the pressure, over a cell of the primal grid, $\Omega_{j}$ for $ j = 1 \ldots n $
\begin{equation}\label{eq:weakINS:p}
	\int_{\Omega_{j}}  \theta\  \nabla \cdot u \  d \textbf{x} = 0,
\end{equation}
where $d \textbf{x} = dxdy$.
\end{subequations}

Integrating by parts the viscous term in \eqref{eq:weakINS:u}, we must take into account that velocity at intercell boundaries is discontinuous and it is necessary to penalize the jumps in order to achieve a stable discretization. 
We associate with this term the bilinear form:
\begin{equation}\label{eq:discr:L}
B(u , \psi) = \int_{\Omega_{i}^{*}}\mu \nabla u  \cdot \nabla \psi \ d \textbf{x} 
+ \, \epsilon \! \int_{\partial \Omega_{i}^{*}}  \mu \lsem u  \rsem \cdot \{\nabla \psi\} \ d\Gamma +\\
- \int_{\partial \Omega_{i}^{*}}  \{\nabla  u \} \cdot \lsem \psi \rsem \ d\Gamma
+ \int_{\partial \Omega_{i}^{*}} \alpha \mu \lsem u  \rsem \lsem \psi \rsem \ d\Gamma,
\end{equation}
where $\alpha = \tfrac{\alpha_{0}}{\dx}$ is the penalization \cite{Arnold:02}.
Changing the sign of $\epsilon$ we obtain symmetric (SIP) \cite{SIP:78} and non-symmetric Interior Penalty (NIP) method \cite{NIP:99}. In the first case the velocity jump term for the mean of the test function is subtracted in the bilinear form, so $\epsilon= - 1$, while in the second method it is added.
Following to \cite{Arnold:02}, the bilinear form $B$ is coercive $\forall \alpha_{0} > 0$ in the NIP case and for $\alpha_{0} > \hat \alpha > 0$, for some $\hat \alpha$ in the SIP case. The estimation of $\hat \alpha$ is in general a nontrivial task, but the advantage of SIP is that the resulting matrix is symmetric and positive definite. Due to the advantage properties of SIP we discretize the viscosity term with this method and for all the test in this article we choose $\alpha_{0} = 1$. 

The integrand of the pressure term in \eqref{eq:weakINS:u} contains
a discontinuity since the pressure is defined on the primal grid and is thus not continuous on the dual velocity cells. The pressure integral is then split
as follows:
\begin{equation}\label{eq:discr:G}
\int_{\Omega_{i}^{*}}  \psi\ \nabla p\  d \textbf{x} =
\int_{\Omega_{i}^{*} \cap \Omega_{i-1} }  \psi\ \nabla p\  d \textbf{x} +
\int_{\Omega_{i}^{*} \cap \Omega_{i}}  \psi\ \nabla p\  d \textbf{x} +
\int_{\Gamma_i }  \psi\  ( p_{i} - p_{i-1} ) d\Gamma
,
\end{equation}
where $p_{i-1}$ and $p_{i}$ denote the discrete pressure in the cells $\Omega_{i-1}$ and $\Omega_{i}$ respectively and $\Gamma_i$ is the interface between $\Omega_{i-1}$ and $\Omega_{i}$, which is located in the middle of $\Omega_{i}^{*}$. 

A similar difficulty appears in \eqref{eq:weakINS:p}, since the discrete velocity is discontinuous on pressure elements, and this is circumvented by computing the divergence term as
\begin{equation}\label{eq:discr:D}
\int_{\Omega_{j}}  \theta\ \nabla\cdot u\  d \textbf{x} =
\int_{\Omega_{j} \cap \Omega_{j}^* }  \theta\ \nabla\cdot u \, d \textbf{x} +
\int_{\Omega_{j} \cap \Omega_{j+1}^*}  \theta\ \nabla\cdot u \, d \textbf{x} +
\int_{\Gamma^*_{j} }  \theta\  ( u_{j+1} - u_{j})  d\Gamma.
\end{equation}
Here above, $\Gamma_j$ denotes the interface between $\Omega_{j}^{*}$ and $\Omega_{j+1}^{*}$, which is located in the middle of $\Omega_{j}$.

Further, for stability, a penalty term must be added to the discretized continuity equation \eqref{eq:weakINS:p} due to the choice of a discontinuous approximation for pressure \cite{Kanschat:07}. Equation \eqref{eq:weakINS:p} is thus modified adding the term
\begin{equation}\label{eq:discr:E}
	\int_{\Gamma_{j}} \alpha  \lsem p \rsem  \lsem \theta \rsem \ d\Gamma
\end{equation}
where the penalization constant is $\alpha = \dx$.
Without this additional term, pressure oscillations that grow as $\dx\to0$ would appear  at the cell interfaces of the main grid.

The left hand side of \eqref{eq:weakINS:u} gives rise to a mass matrix term and to a convective term that depends nonlinearly on $u$. By considering in \eqref{eq:weakINS} an implicit discretization for all terms except for the nonlinear convective term, one obtains a linear system for the velocity and pressure unknowns at time $t^{n+1}$ that has the following block structure

\begin{equation}\label{eq:coeffm}
\mathcal{A} \textbf{x} = \textbf{f} \iff
\begin{bmatrix}
N & G \\
D & E \\
\end{bmatrix}
\begin{pmatrix}
u \\
p \\
\end{pmatrix}
= 
\begin{pmatrix}
b_u(u) \\
0 \\
\end{pmatrix}
.
\end{equation}
Here above, $ N = M + L$ is a square matrix formed by  $L$ and $M$ that discretize the Laplacian and the mass operator; these are of size $\Ogrande(1)$ and $\Ogrande(\dx)$ respectively.
$G$ is a rectangular tall matrix of size $\Ogrande(\dt)$ corresponding to the gradient operator \eqref{eq:discr:G}, while $D$, coming from \eqref{eq:discr:D}, is its transpose up to a scaling factor, which has size $\Ogrande(1)$.
Finally $E$ is a square matrix of size $\Ogrande(\dx)$ containing the penalty term \eqref{eq:discr:E}.
In the right hand side, $b_u(u)$ is the discretization of the nonlinear convective terms with a classical explicit TVD Runge–Kutta method and Rusanov fluxes, as in \cite{Tavelli:14}. Boundary conditions for a prescribed velocity profile at the inlet are inserted in the system in place of the first rows of $N, G$ and $b(u)$; we impose an outlet pressure by prescribing the stress modifying the last rows of the same blocks.

The time step $\dt$ is restricted by a CFL-type restriction for DG schemes depending only on the fluid velocity.
In the following analysis, we thus assume that $\frac{\dt}{\dx}=c=\Ogrande(1)$.


\section{Preliminaries}

\label{sec:prelim}

Here we first formalize the definition of block Toeplitz and circulant sequences associated to a matrix-valued Lebesgue integrable function (see Subsection \ref{sub:toeplitz}). Moreover, in Subsection \ref{sub:glt} we introduce a class of matrix-sequences containing block Toeplitz sequences known as the block Generalized Locally Toeplitz (GLT) class \cite{GS,axioms,barbarino2020block}. The properties of block GLT sequences and few other new spectral tools introduced in Subsection \ref{sub:new-tools} will be used to derive the spectral properties  of $\mathcal{A}$ in \eqref{eq:coeffm} as well as of its blocks and its Schur complement. 

\subsection{Block Toeplitz and circulant matrices}\label{sub:toeplitz}

Let us denote by $L^1([-\pi,\pi],s)$ the space of $s\times s$ matrix-valued functions $f:[-\pi,\pi]\rightarrow\mathbb{C}^{s\times s}$, $f=[f_{ij}]_{i,j=1}^s$ with $f_{ij}\in L^1([-\pi,\pi])$, $i,j=1,\dots,s$. In Definition \ref{def:Tblock} we introduce the notion of Toeplitz and circulant matrix-sequences generated by $f$.
\begin{definition}\label{def:Tblock}
	Let $f\in L^1([-\pi,\pi],s)$ and let $t_j$ be its Fourier coefficients
	\begin{equation*}
	t_{j}:=\frac1{2\pi}\int_{-\pi}^\pi f(\theta){\rm e}^{-{\bf i}j\theta}\ \dd \theta\in\mathbb{C}^{s\times s},
	\end{equation*}
	where the integrals are computed component-wise. 
	Then, the $n$-th $s\times s$-block Toeplitz matrix associated with $f$ is the matrix of order $\widehat n=s\cdot n$ given by
	\begin{equation*}
	T_{n}(f)=\left[t_{i-k}\right]_{i,k=1}^{n}.
	\end{equation*}
	Similarly, the $n$-th $s\times s$-block circulant matrix associated with $f$ is the following $\widehat n\times \widehat n$ matrix
	\begin{align*}
	C_{n}(f)=\left[t_{(i-k){\rm mod}n}\right]_{i,k=1}^{n}.
	\end{align*}
	The sets $\{T_{n}(f)\}_{{ n}}$ and $\{C_{n}(f)\}_{{ n}}$ are
	called the \emph{families of $s\times s$-block Toeplitz and circulant matrices generated by $f$}, respectively. The function $f$ is referred to as the \emph{generating function} either of $\{T_{n}(f)\}_{{n}}$ or $\{C_{n}(f)\}_{{n}}$.
\end{definition}

It is useful for our later studies to extend the definition of block-Toeplitz sequence also to the case where the symbol is a rectangular matrix-valued function.
\begin{definition}\label{def:Tblock-rect}
Let $f:[-\pi,\pi]\rightarrow\C^{s\times q}$, with $s\ne q$, and such that $f_{ij}\in L^1([-\pi,\pi])$ for $i=1,\ldots,s$ and $j=1,\ldots,q$. Then, given $n\in\mathbb{N}$, we denote by $T_{n}(f)$ the $s\cdot n\times q\cdot n$ matrix whose entries are $T_{n}(f)=[t_{i-k}]_{i,k=1}^{n}$, with $t_{j}\in\C^{s\times q}$ the Fourier coefficients of $f$.	
\end{definition}

The generating function $f$ provides a description of the spectrum of $T_{n}(f)$, for $n$ large enough in the sense of the following definition.

\begin{definition}\label{def-distribution}
	Let $f:[a,b]\to\mathbb{C}^{s\times s}$ be a measurable matrix-valued function with eigenvalues $\lambda_i(f)$ and singular values $\sigma_i(f)$, $i=1,\ldots,s$. Assume that $\{A_n\}\p$ is a sequence of matrices such that ${\rm dim}(A_n)=d_n\rightarrow\infty$, as $n\rightarrow\infty$ and with eigenvalues $\lambda_j(A_n)$ and singular values $\sigma_j(A_n)$, $j=1,\ldots,d_n$. 
	\begin{itemize}
		\item We say that $\{A_n\}\p$ is {\em distributed as $f$ over $[a,b]$ in the sense of the eigenvalues,} and we write $\{A_n\}\p\sim_\lambda(f,[a,b]),$ if
		\begin{equation}\label{distribution:eig}
		\lim_{n\to\infty}\frac{1}{d_n}\sum_{j=1}^{d_n}F(\lambda_j(A_n))=
		\frac1{b - a} \int_a^b \frac{\sum_{i=1}^sF(\lambda_i(f(t)))}{s} \dd t,
		\end{equation}
		for every continuous function $F$ with compact support. In this case, we say that $f$ is the \emph{spectral symbol} of $\{A_{n}\}_{n}$.
		
		\item We say that $\{A_n\}\p$ is {\em distributed as $f$ over $[a,b]$ in the sense of the singular values,} and we write $\{A_n\}\p\sim_\sigma(f,[a,b]),$ if
		\begin{equation}\label{distribution:sv}
		\lim_{n\to\infty}\frac{1}{d_n}\sum_{j=1}^{d_n}F(\sigma_j(A_n))=
		\frac1{b - a} \int_a^b \frac{\sum_{i=1}^sF(\sigma_i(f(t)))}{s} \dd t,
		\end{equation}
		for every continuous function $F$ with compact support.
	\end{itemize}	
\end{definition}

Throughout the paper, when the domain can be easily inferred from the context, we replace the notation $\{A_n\}_n\sim_{\lambda,\sigma}(f,[a,b])$ with $\{A_n\}_n\sim_{\lambda,\sigma} f$.

\begin{remark}\label{rem:approx}
	If $f$ is smooth enough, an informal interpretation of the limit relation
	\eqref{distribution:eig} (resp. \eqref{distribution:sv})
	is that when $n$ is sufficiently large,
	then $d_n/s$ eigenvalues (resp. singular values) of $A_{n}$ can
	be approximated by a sampling of $\lambda_1(f)$ (resp. $\sigma_1(f)$)
	on a uniform equispaced grid of the domain $[a,b]$, and so on until the
	last $d_n/s$ eigenvalues (resp. singular values), which can be approximated by an equispaced sampling
	of $\lambda_s(f)$ (resp. $\sigma_s(f)$) in the domain.
\end{remark}

For Toeplitz matrix-sequences, the following theorem due to Tilli  holds, which generalizes previous researches along the last 100 years by Szeg\H{o}, Widom, Avram, Parter, Tyrtyshnikov, Zamarashkin (see \cite{barbarino2020block,bottcher2013analysis,GS,tyrtyshnikov1998spectra} and references therein).

\begin{theorem}[see \cite{Tillinota}]\label{szego-herm}
	Let $f\in L^1([-\pi,\pi],s)$, then $\{T_{n}(f)\}_{{n}}\sim_\sigma(f,[-\pi, \pi]).$ If $f$ is a Hermitian matrix-valued function, then $\{T_{n}(f)\}_{{n}}\sim_\lambda(f,[-\pi, \pi]).$
\end{theorem}

Since rectangular matrices always admit a singular value decomposition, equation \eqref{distribution:sv} can also be extended to rectangular matrix-sequences. Throughout we denote by $A_{m_1,m_2,s,q}\in\C^{s\cdot m_1\times q \cdot m_2}$ the rectangular matrix that has $m_1$ blocks of $s$ rows and $m_2$ blocks of $q$ columns. As a special case, with $[T_n(f)]_{m_1,m_2,s,q}$, $m_1,m_2\le n$ we denote the \lq leading principal' submatrix of $T_n(f)$ of size $s\cdot m_1\times q \cdot m_2$. Moreover, if $f\in\C^{s\times q}$ then we omit the subscripts $s,q$ since they are implicitly clear from the size of the symbol.

\begin{definition}\label{def:distr-rect}
Given a measurable function $f:[a,b]\rightarrow\C^{s\times q}$, with $s\ne q$ and a matrix-sequence $\{A_{m_1,m_2,s,q}\}_n$, with $A_n\in\C^{s \cdot m_1\times q\cdot m_2 }$, $m_1\sim m_2$, $m_1,m_2\rightarrow\infty$ as $n\rightarrow\infty$ then we say that $\{A_{m_1,m_2,s,q}\}_n\sim_{\sigma} (f,[a,b])$ iff
	\begin{equation*}
	\lim_{n\to\infty}\frac{1}{s\cdot m_1\land q\cdot m_2}\sum_{j=1}^{s\cdot m_1\land q\cdot m_2}F(\sigma_j(A_{m_1,m_2,s,q}))=
	\frac1{b - a} \int_a^b \frac{\sum_{i=1}^{s\land q}F(\sigma_i(f(t)))}{s\land q} \dd t,
	\end{equation*}
	with $x\land y:=\min\{x,y\}$, for every continuous function $F$ with compact support.
\end{definition}

\begin{remark} \label{rem:distr-rect}
Based on Definition \ref{def:distr-rect} the first part of Theorem \ref{szego-herm} extends also to rectangular block Toeplitz matrices in the sense of Definition \ref{def:Tblock-rect} (see \cite{Tillinota}) as well as to sequences whose $n$-th matrix is $A_{m_1,m_2,s,q}=[T_n(f)]_{m_1,m_2}$, $f\in\C^{s\times q}$, with $m_1,m_2\le n$, $m_1\sim m_2$, $m_1,m_2\rightarrow \infty$ as $n\rightarrow\infty$.
\end{remark}

The following theorem is a useful tool for computing the spectral distribution of a sequence of Hermitian matrices. For the related proof, see \cite[Theorem 4.3]{mari}. Here, the conjugate transpose of the matrix $X$ is denoted by $X^*$.

\begin{theorem}\label{th:extradimensional}
	Let $\{A_n\}_n$ be a sequence of matrices, with $A_n$ Hermitian of size $d_n$, and let $\{P_n\}_n$ be a sequence such that $P_n\in\mathbb C^{d_n\times\delta_n}$, $P_n^*P_n=I_{\delta_n}$, $\delta_n\le d_n$ and $\delta_n/d_n\to1$ as $n\to\infty$. Then $\{A_n\}_n\sim_{\lambda}f$ if and only if $\{P_n^*A_nP_n\}_n\sim_{\lambda}f$.
\end{theorem}

The following result allows us to determine the spectral distribution of a Hermitian matrix-sequence plus a correction (see \cite{barbarino2018non}).

\begin{theorem}\label{thm:quasi-herm}
	Let $\{X_n\}_n$ and $\{Y_n\}_n$ be two matrix-sequences, with $X_n,Y_n\in\C^{d_n\times d_n}$, and assume that 
	\begin{itemize}
		\item[{(a)}] $X_n$ is Hermitian for all $n$ and $\{X_n\}_n\sim_{\lambda}f$;
		\item[{(b)}] $\|Y_n\|_F=o(\sqrt{d_n})$ as $n\rightarrow\infty$, with $\|\cdot\|_{F}$ the Frobenius norm.
	\end{itemize}
	Then, $\{X_n+Y_n\}_n\sim_{\lambda}f$.
\end{theorem}

For a given matrix $X\in\C^{m\times m}$, let us denote by $\|X\|_{1}$ the trace norm defined by $\|X\|_{1}:=\sum_{j=1}^{m}\sigma_j(X)$, where $\sigma_j(X)$ are the $m$ singular values of $X$.

\begin{corollary}\label{cor:quasi-herm}
	Let $\{X_n\}_n$ and $\{Y_n\}_n$ be two matrix-sequences, with $X_n,Y_n\in\C^{d_n\times d_n}$, and assume that (a) in Theorem~\ref{thm:quasi-herm} is satisfied. Moreover, assume that any of the following two conditions is met:
	\begin{itemize}
		\item $\|Y_n\|_{1}=o(\sqrt{d_n})$;
		\item $\|Y_n\|=o(1)$, with $\|\cdot\|$ being the spectral norm.
	\end{itemize}
	Then, $\{X_n+Y_n\}_n\sim_{\lambda}f$.	
\end{corollary}

We end this subsection by reporting the key features of the block circulant matrices, also in connection with the generating function.

\begin{theorem}[\cite{garoniserrasesana15}]\label{circ-s}
	Let $f\in L^1([-\pi,\pi],s)$ be a matrix-valued function with $s\ge 1$ and let $\{t_j\}_{j\in\mathbb{Z}}$, $t_j\in\mathbb{C}^{s\times s}$ be its Fourier coefficients. Then, the following (block-Schur) decomposition of 
	$C_{n}(f)$ holds:
	\begin{equation}
	\label{schur-s}
	C_{n}(f)=(F_{n}\otimes I_s) D_{n}(f) (F_{n}\otimes I_s)^*, 
	\end{equation}
	where 	
		\begin{equation}\label{eig-circ}
	D_{n}(f) =\diag_{0\le r\le  n-1}\left(S_{n}(f)\left(\theta_{r}\right)\right), \quad \theta_{r}=\frac{2\pi r}{n}, \quad F_{n}=\frac{1}{\sqrt{n}} \left({\rm e}^{-{\ii}j\theta_{r}}\right)_{j,r=0}^{n-1}
	\end{equation}
	with $S_{n}(f)(\cdot)$ the $n$-th Fourier sum of $f$ given by
	\begin{equation}\label{fourier-sum}
	S_{n}(f)(\theta) = \sum_{j=0}^{n-1} t_{j} 
	{\rm e}^{\ii j\theta}.
	\end{equation}
	Moreover, the eigenvalues of $C_{n}(f)$ are given by the evaluations of $\lambda_t(S_{n}(f)(\theta))$, $t=1,\ldots,s$, if $s\ge2$ or of $S_{n}(f)(\theta)$ if $s=1$ at the grid points $\theta_r$.
\end{theorem}

\begin{remark}\label{fourier-vs-funzione}
	If $f$ is a trigonometric polynomial of fixed degree (with respect to $n$), then it is worth noticing that
	$S_{n}(f)(\cdot) = f(\cdot)$ for $n$ large enough: more precisely, $n$ should be larger than the double of the degree.
	Therefore, in such a setting, the eigenvalues of $C_{n}(f)$ are either the evaluations of $f$ at the grid points if $s=1$ or
	the evaluations of $\lambda_t(f(\cdot))$, $t=1,\ldots,s$, at the very same grid points.
\end{remark}

We recall that every matrix/vector operation with circulant matrices has cost $O(\widehat{n}\log \widehat {n})$ with moderate multiplicative constants: in particular, this is true for the matrix-vector product, for the solution of a linear system, for the computation of the blocks
$S_{n}(f)(\theta_{r})$ and consequently of
the eigenvalues (see e.g. \cite{fft}).

\subsection{Block Generalized locally Toeplitz class}\label{sub:glt}
\noindent
In the sequel, we introduce the block GLT class, a $*$-algebra of matrix-sequences containing block Toeplitz matrix-sequences. The formal definition of block GLT matrix-sequences is rather technical, therefore we just give and briefly discuss a few properties of the block GLT class, which are sufficient for studying the spectral features of $\mathcal{A}$ as well as of its blocks and its Schur complement.

Throughout, we use the following notation
$$\{A_n\}\p\sim_{\rm GLT} { \kappa(x,\theta)},\quad \ \ \kappa:[0,1]\times[-\pi,\pi]\rightarrow\mathbb{C}^{s\times s},$$
to say that the sequence $\{A_n\}\p$ is a $s\times s$-block GLT sequence with GLT symbol $\kappa(x,\theta)$.

Here we list four main features of block GLT sequences.
\begin{itemize}
	\item[{\bf GLT1}] Let $\{A_{ n}\}_{ n}\sim_{\rm GLT}\kappa$ with $\kappa:G\rightarrow \mathbb{C}^{s\times s}$, $G=[0,1]\times[-\pi,\pi]$, then $\{A_{ n}\}_{ n}\sim_\sigma(\kappa,G)$. If the matrices $A_{ n}$ are Hermitian, then it also holds that $\{A_{ n}\}_{ n}\sim_\lambda(\kappa,G)$.
	\item[{\bf GLT2}] The set of block GLT sequences forms a $*$-algebra, i.e., it is closed under linear combinations, products, conjugation, but also inversion when the symbol is invertible a.e. In formulae, let $\{ A_{ n} \}_{ n} \glteq \kappa_1$ and $\{ B_{ n} \}_{ n} \glteq \kappa_2$, then
	\begin{itemize}
		
		\item[$\bullet$] $\{\alpha A_{ n} + \beta B_{ n}\}_{ n} \glteq \alpha\kappa_1+\beta \kappa_2, \quad \alpha, \beta \in \mathbb{C};$
		\item[$\bullet$] $\{A_{ n}B_{ n}\}_{ n} \glteq \kappa_1 \kappa_2;$
		\item[$\bullet$] $\{ A_{ n}^{*} \}_{ n} \glteq {\kappa^*_1};$
		\item[$\bullet$] $\{A^{-1}_{ n}\}_{ n} \glteq \kappa_1^{-1}$ provided that $\kappa_1$ is invertible a.e.
	\end{itemize}
	\item[{\bf GLT 3}] Any sequence of block Toeplitz matrices $\{ T_{ n}(f) \}_{ n}$ generated by a function $f \in  L^1([-\pi, \pi],s)$ is a $s\times s$-block GLT sequence with symbol $\kappa(x, \theta) = f(\theta)$.
	\item[{\bf GLT4}] Let $\{A_{ n}\}_{ n}\sim_\sigma 0$. We say that $\{A_{ n}\}_{ n}$ is a \emph{zero-distributed matrix-sequence}. Note that for any $s>1$ $\{A_{ n}\}_{ n}\sim_\sigma O_s$, with $O_s$ the $s\times s$ null matrix, is equivalent to $\{A_{ n}\}_{ n}\sim_\sigma 0$. Every zero-distributed matrix-sequence is a block GLT sequence with symbol $O_s$ and viceversa, i.e., $\{A_{ n}\}_{ n}\sim_\sigma0$ $\iff$ $\{A_{ n}\}_{ n}\sim_{\rm GLT}O_s$.
\end{itemize}

According to Definition \ref{def-distribution}, in the presence of a zero-distributed sequence 
the singular values of the $n$-th matrix (weakly) cluster around $0$. This is formalized in the following result \cite{GS}.

\begin{proposition}\label{prop:zero}
	Let $\{ A_n \}_n$ be a matrix sequence with $A_n$ of size $d_n$ with $d_n\rightarrow\infty$, as $n\rightarrow\infty$. Then
	$\{ A_n \}_n \sim_{\sigma} 0$ if and only if there
	exist
	two matrix sequences $\{ R_n \}_n$ and $\{ E_n\}_n$ such that
	$A_n = R_n + E_n$, and
	\[
	\lim_{n \to \infty} \frac{\mathrm{rank}(R_n)}{d_n} = 0, \qquad
	\lim_{n \to \infty} \|E_n\| = 0.
	\]
	The matrix $R_n$ is called rank-correction and the matrix $E_n$ is called norm-correction.
\end{proposition}

\subsection{Some new spectral tools}\label{sub:new-tools}

In this subsection we introduce some new spectral tools that will be used in Section \ref{sec:spectr}.

The following theorem concerns the spectral behavior of matrix-sequences whose $n$-th matrix is a product of a square block Toeplitz matrix by a rectangular one. 

\begin{theorem}\label{th:distr-rect}
Let $f:[-\pi,\pi]\rightarrow\C^{s\times s}$ and let $g:[-\pi,\pi]\rightarrow\C^{s\times q}$, $h:[-\pi,\pi]\rightarrow\C^{q\times s}$ with $q< s$. Then 
\begin{equation}\label{eq:distr-rect1}
\{T_n(f)T_{n}(g)\}_n\sim_{\sigma}(f\cdot g,[-\pi,\pi]),
\end{equation} 
and 
\begin{equation}\label{eq:distr-rect2}
\{T_n(h)T_{n}(f)\}_n\sim_{\sigma}(h\cdot f,[-\pi,\pi]).
\end{equation}
\end{theorem}
{\bf Proof.} We only prove relation \eqref{eq:distr-rect1}, since the same argument easily brings to \eqref{eq:distr-rect2} as well. Let us define $g_{\rm ex}:[-\pi,\pi]\rightarrow\C^{s\times s}$ obtained completing $g$ with $s-q$ null columns. By \textbf{GLT3} and \textbf{GLT2} we know that
\begin{equation}\label{eq:distr_ex}
\{T_n(f)T_{n}(g_{\rm ex})\}_n\sim_{\sigma}(f\cdot g_{\rm ex},[-\pi,\pi]).
\end{equation}
Let us now explicitly write \eqref{eq:distr_ex} according to Definition \ref{def-distribution}
\begin{align*}
 \lim_{n\to\infty}\frac{1}{sn}\sum_{j=1}^{sn}F(\sigma_j(T_n(f)T_{n}(g_{\rm ex})))&=
\frac1{2\pi} \int_{-\pi}^{\pi} \frac{\sum_{i=1}^sF(\sigma_i(f(t)g_{\rm ex}(t)))}{s} \dd t.
\end{align*}
The left-hand side of the previous equation can be rewritten as follows
\begin{align*}
 \lim_{n\to\infty}\frac{1}{sn}\sum_{j=1}^{sn}F(\sigma_j(T_n(f)T_{n}(g_{\rm ex})))
 &=\lim_{n\to\infty}\frac{1}{sn}\left[\sum_{j=1}^{qn}F(\sigma_j(T_n(f)T_{n}(g_{\rm ex})))+\sum_{qn+1}^{sn}F(0)\right]\\
 &=\lim_{n\to\infty}\frac{1}{sn}\sum_{j=1}^{qn}F(\sigma_j(T_n(f)T_{n}(g)))+ \frac{(s-q)}{s}F(0),
\end{align*}
while manipulating the right-hand side we obtain
\begin{align*}
\frac1{2\pi} \int_{-\pi}^{\pi} \frac{\sum_{i=1}^sF(\sigma_i(f(t)g_{\rm ex}(t)))}{s} \dd t&=
\frac1{2\pi} \int_{-\pi}^{\pi} \frac{\sum_{i=1}^qF(\sigma_i(f(t)g_{\rm ex}(t))) +\sum_{i=q+1}^sF(0)}{s} \dd t\\
&=\frac1{2\pi} \int_{-\pi}^{\pi} \frac{\sum_{i=1}^qF(\sigma_i(f(t)g(t))) +(s-q)F(0)}{s} \dd t\\
&=\frac{1}{2\pi} \int_{-\pi}^{\pi} \frac{\sum_{i=1}^qF(\sigma_i(f(t)g(t)))}{s}\dd t +\frac{(s-q)}{s}F(0).
\end{align*}
Therefore we arrive at 
\begin{align*}
\lim_{n\to\infty}\frac{1}{sn}\sum_{j=1}^{qn}F(\sigma_j(T_n(f)T_{n}(g)))&=
\frac1{2\pi} \int_{-\pi}^{\pi} \frac{\sum_{i=1}^qF(\sigma_i(f(t)g(t)))}{s} \dd t.
\end{align*}
which proves \eqref{eq:distr-rect1}, once multiplied by $\frac{s}{q}$.\hfill $\square$

\begin{remark}\label{rem:th-ext}
Theorem \ref{th:distr-rect} can easily be extended to the case where also $T_n(f)$ is a properly sized rectangular block Toeplitz matrix. In particular, when 
$f\cdot g$ (or $h \cdot f$) results in a Hermitian square matrix-valued function then the distribution also holds in the sense of the eigenvalues.
\end{remark}

Along the same lines of the previous theorem the following result holds. We notice that Theorem \ref{th:distr-rect} and Theorem \ref{th:distr-rect-bis} are special cases of a more general theory which connects GLT sequences having symbols with different matrix sizes (see \cite{GLTheretogeneous}). 

\begin{theorem}\label{th:distr-rect-bis}
Let $g:[-\pi,\pi]\rightarrow\C^{s\times s}$ be Hermitian positive definite almost everywhere and let $f:[-\pi,\pi]\rightarrow\C^{q\times s}$ with $q< s$. Then 
\begin{equation*}
\{T_n(f)T_{n}^{-1}(g)T_n(f^*)\}_n\sim_{\sigma}(f\cdot g^{-1}\cdot f^*,[-\pi,\pi]),
\end{equation*} 
and 
\begin{equation*}
\{T_n(f)T_{n}^{-1}(g)T_n(f^*)\}_n\sim_{\lambda}(f\cdot g^{-1}\cdot f^*,[-\pi,\pi]).
\end{equation*}
\end{theorem}

The following theorem will be used in combination with Theorem \ref{th:extradimensional} to obtain the spectral symbol of the whole coefficient matrix sequence appearing in \eqref{eq:coeffm}. 
\begin{theorem} \label{th:wholeA}
Let 
\begin{equation*}
A_n=
\begin{bmatrix}
T_{n}(\mathfrak{f}_{11}) & T_{n}(\mathfrak{f}_{12})\\
T_{n}(\mathfrak{f}_{21}) & T_{n}(\mathfrak{f}_{22})
\end{bmatrix}
\end{equation*}
with $\mathfrak{f}_{11}:[-\pi,\pi]\rightarrow \C^{k\times k}$, $\mathfrak{f}_{12}:[-\pi,\pi]\rightarrow \C^{k\times q}$, $\mathfrak{f}_{21}:[-\pi,\pi]\rightarrow \C^{q\times k}$, $\mathfrak{f}_{22}:[-\pi,\pi]\rightarrow \C^{q\times q}$, $k,q\in\mathbb{N}$. Then there exists a permutation matrix $\Pi$ such that $A_n=\Pi T_n(\mathfrak{f}) \Pi^T$ with 
\begin{equation*}
\mathfrak{f}=
\begin{bmatrix}
	\mathfrak{f}_{11} & \mathfrak{f}_{12}\\
	\mathfrak{f}_{21} & \mathfrak{f}_{22}
\end{bmatrix}.
\end{equation*}
Hence $A_n$ and $T_n(\mathfrak{f})$ share the same eigenvalues and the same singular values and consequently $\{A_n\}_n$ and $\{T_n(\mathfrak{f})\}_n$ enjoy the same distribution features.
\end{theorem}

\textbf{Proof.} Let $I_{kn+qn}$ be the identity matrix of size $kn+qn$ and let us define the following sets of indexes $H=\{1,\ldots,kn+qn\}$ and $J=\{k+1,\ldots,k+q,2k+q+1,\ldots,2k+2q,3k+2q+1,\ldots,3k+3q,\ldots,nk+(n-1)q+1,\ldots,nk+nq\}$. Let $\Pi$ be the $(kn+qn)\times(kn+qn)$-matrix whose first $kn$ rows are defined as the rows of $I_{kn+qn}$ that correspond to the indexes in $H\backslash J$ and the remaining as the rows of $I_{kn+qn}$ that correspond to the indexes in $J$. The thesis easily follows observing that $\Pi$ is the permutation matrix that relates $A_n$ and $T_n(\mathfrak{f})$. 

Thus $A_n$ and $T_n(\mathfrak{f})$ are similar because $\Pi^T$ is the inverse of $\Pi$ and as consequence both matrices $A_n$ and $T_n(\mathfrak{f})$ share the same eigenvalues. Furthermore both $\Pi$ and $\Pi^T$ are unitary and consequently by the singular value decomposition the two  matrices $A_n$ and $T_n(\mathfrak{f})$ share the same singular values. Finally it is transparent that one of the matrix sequences (between  $\{A_n\}_n$ and $\{T_n(\mathfrak{f})\}_n$) has a distribution if and only the other has the very same distribution.
$\square$

\section{Spectral analysis}\label{sec:spectr}

This section concerns the spectral study of the matrix $\mathcal{A}$ in \eqref{eq:coeffm} together with its blocks and Schur complement. In the following, we consider the case of $d(x)=d$ (constant width); we choose at first the smallest nontrivial case which is $n_x=1$ and $n_y=3$ ($n_u=(n_x+1)(n_y-1)=4$ and $n_p=(n_x+1)=2$) and then comment on the general case.

\subsection{Spectral study of the blocks of $\mathcal{A}$}\label{sub:blocks}
\noindent We start by spectrally analyzing the four blocks that compose the matrix $\mathcal{A}$.

\paragraph{Laplacian and mass operator}
The $(1,1)$ block $N$ of $\mathcal{A}$ in \eqref{eq:coeffm} is a sum of two terms: the Laplacian matrix $L$ and the mass matrix $M$ that are respectively obtained by testing the PDE term $\nabla \cdot(\mu \ \nabla u)$ and the term $\partial_{t} u$ with the basis functions for velocity.
  	
\begin{figure}[!h]
	\begin{center}
		\includegraphics{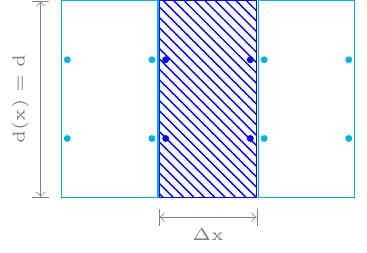}
	\end{center}
	\caption{Illustration of the stencil that refers to the mass and Laplacian matrix.}
	\label{fig:discretL}
\end{figure}

The matrix $L$
is organized in blocks of rows each of size $n_u=4$ which corresponds to the number of test functions per cell (associated with the blue degrees of freedom in Fig.~\ref{fig:discretL}); in each row there are at most twelve nonzeros elements 
(associated with all the degrees of freedom in Fig.~\ref{fig:discretL}). Using SIP in \eqref{eq:discr:L} and excluding the boundary conditions, we can write
\[
L_{n+1} = \frac{27}{70}d \mu c U_{n+1}
\]
with
\begin{eqnarray*}
	U_{n+1} & = & 
	\tridiag
	\left[
		\begin{array}{cccc|cccc|cccc}\notag
		-\tfrac{1}{ 2}	& \tfrac{ 1}{16} &  0 & 0 &  1 & -\tfrac{1}{8} & 0 & 0 &  -\tfrac{1}{ 2}  & \tfrac{1}{16} & 0 & 0 \\
		\tfrac{ 1}{16}	& -\tfrac{1}{ 2} &  0 & 0 &  -\tfrac{1}{8} & 1 & 0 & 0 &   \tfrac{1}{16}	& -\tfrac{1}{2} & 0 & 0 \\
		0 & 0 &  -\tfrac{1}{ 2}	& \tfrac{ 1}{16} &  0 & 0 & 1 & -\tfrac{1}{8} &   0 & 0 & -\tfrac{1}{ 2}	& \tfrac{ 1}{16} \\
		0 & 0 &  \tfrac{1}{ 16}	&  -\tfrac{ 1}{2} &  0 & 0 & -\tfrac{1}{8} & 1 &   0 & 0 &  \tfrac{ 1}{16}& -\tfrac{1}{ 2}
		\end{array}
	\right] \\
		& & +\Ogrande(\dx^2),
\end{eqnarray*}
where $\mu$ is the viscosity, $c = \frac{\dt}{\dx}$, and $n+1$ is the number of velocity cells.

It is then clear that $L_{n+1}$
is a $4\times4$-block Toeplitz matrix of size $\widehat{n}=4 \cdot (n+1)$.
As a consequence, we can obtain insights on its spectrum studying the symbol associated to $\{L_{n+1}\}_{n}$. With this aim, let us define
\begin{equation}\label{eq:X}
	X = 
	\begin{bmatrix}\notag
	\tfrac{ 1}{ 2}  &  -\tfrac{ 1}{16} \\
	-\tfrac{ 1}{16} &    \tfrac{ 1}{ 2} \\
	\end{bmatrix}
	,
\end{equation}
and $l_1,l_0,l_{-1}$ as follows
\[
	l_{1} =
	\left[
		\begin{array}{c|c}\notag
		-X  &  0 \\
		\hline
		0  &   -X \\
		\end{array}
	\right],
	\hspace{20pt}
	l_{0} =
	\left[
		\begin{array}{c|c}\notag
		2X  &  0 \\  
		\hline
		0  &   2X \\
		\end{array}
	\right],
	\hspace{20pt}
	l_{-1} =
	\left[
		\begin{array}{c|c}\notag
		-X  &  0 \\
		\hline
		-0  &   -X \\
		\end{array}
	\right].
\]
Since we are assuming that $c=\Ogrande(1)$ the symbol associated to $\{L_{n+1}\}_n$ is the function $\lsym:[-\pi ,\pi] \to \mathbb{C}^{4\times4}$ defined as
\[
	\lsym(\theta)  = \frac{27}{70} d \mu c (l_{0} + l_{1} e^{\ii\theta} + l_{-1} e^{-\ii\theta})=		\frac{27}{70} d \mu c
	\begin{bmatrix}
	(2-2\cos\theta) & 0 \\
	0 & (2-2\cos\theta) \\
	\end{bmatrix}
	\otimes X.
\]
Recalling Theorem \ref{szego-herm} and \textbf{GLT3}, we conclude that
\begin{equation}\label{eq:distrL}
	\{L_{n+1}\}_n
	\sim_{{\rm GLT},\sigma,\lambda}
	(\lsym,[-\pi,\pi]).
\end{equation}
	
\begin{remark}\label{rem:BCs}
We have assumed that $L_{n+1}$ does not contain the boundary conditions, but if we let them come into play, then the spectral distribution would remain unchanged. Indeed, the matrix that corresponds to the Laplacian operator can be expressed as the sum $L_{n+1}+R_{n+1}$ with  $R_{n+1}$ a rank-correction. Since the boundary conditions imply a correction in a constant number of entries and since the absolute values of such corrections are uniformly bounded with respect to the matrix size, it easily follows that $\|R_{n+1}\|= \Ogrande(1)$ and hence Theorem \ref{thm:quasi-herm} can be applied.
\end{remark}

\begin{figure}[!t]
	\begin{center}
		\begin{tabular}{cc}
			\includegraphics[width=0.4\textwidth]{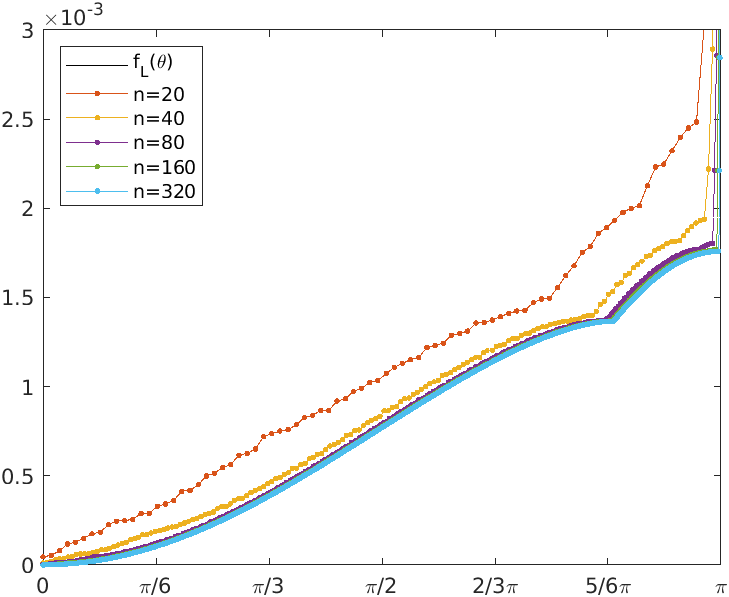}&
			\includegraphics[width=0.4\textwidth]{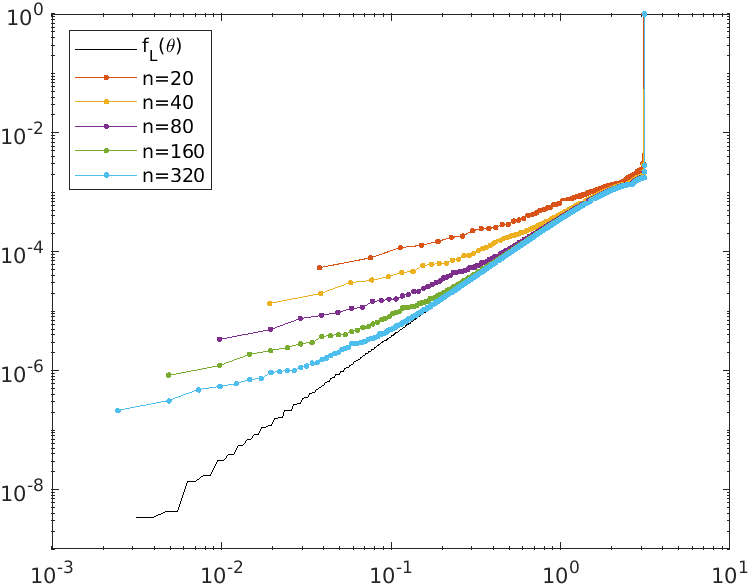}\\
			(a) & (b)
		\end{tabular}
	\end{center}
	\caption{(a) The spectrum of $L_{n+1}$  with different number of cells vs sampling of the eigenvalue functions of the symbol $\lsym(\theta)$; (b) is the same picture, but in bilogarithmic scale.}
	\label{fig:symbolL}
\end{figure}

It is easy to compute the four eigenvalue functions of $\lsym(\theta)$, which are  $\frac{27}{70} d \mu c  2(1-\cos\theta)\left(\tfrac{1}{2}\pm\tfrac{1}{16}\right)$, each with multiplicity 2. Note that all eigenvalue functions vanish at $\theta = 0$ with a zero of second order. Recalling Remark \ref{rem:approx}, we expect that a sampling of the eigenvalues of $\lsym(\theta)$ provides an approximation of the spectrum of the discretized Laplacian operator. This is confirmed in Fig.~\ref{fig:symbolL}, where we compare the Laplacian matrix, including the boundary conditions, with an equispaced sampling of the eigenvalue functions of $\lsym(\theta)$ in $[-\pi,\pi]$.

\begin{figure}[!t]
	\begin{center}
		\begin{tabular}{c}
			\includegraphics[width=0.4\textwidth]{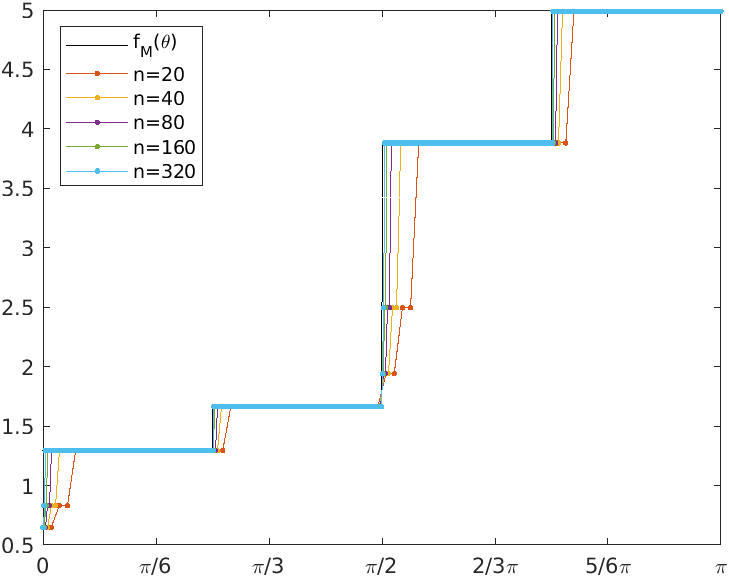}
		\end{tabular}
	\end{center}
	\caption{The eigenvalues of $\frac{1}{\dx} M_{n+1}$ matrix with different number of cells vs sampling of the eigenvalue functions of $\msym(\theta) $.} 
	\label{fig:symbolM}
\end{figure}

The mass matrix $M_{n+1}$ is block diagonal and  has the form
\[
	M_{n+1} = \dfrac{9}{70} d \dx \rho
	\diag
	\left[
	\begin{array}{cccc}\notag
	1	           & -\tfrac{1}{8} &  \tfrac{1}{2} & -\tfrac{1}{16}\\
	-\tfrac{1}{8} & 1             & -\tfrac{1}{16}&  \tfrac{1}{2} \\
	\tfrac{1}{2} & -\tfrac{1}{16}&   1 		   & -\tfrac{1}{8} \\
	-\tfrac{1}{16}&  \tfrac{1}{2} & -\tfrac{1}{8} & 			 1 \\
	\end{array}
	\right].
\]
As for $L_{n+1}$, also $M_{n+1}$ is a  $4\times 4$-block Toeplitz of size $\widehat{n}=4\cdot (n+1)$. In order to study its symbol we look at the scaled matrix-sequence $\{\frac{1}{\Delta x} M_{n+1}\}_n$. The reason for such scaling is that the symbol is defined for sequences of Toeplitz matrices whose elements do not vary with their size. The symbol of the scaled mass-matrix sequence $\{\frac{1}{\Delta x}M_{n+1}\}_n$ can be written as
\[
	\msym(\theta)=\dfrac{9}{70} d \rho
	\left[
	\begin{array}{cc}
	2 & 1 \\
	1 & 2 \\
	\end{array}
	\right]
	\otimes X
\]
with $X$ as in \eqref{eq:X} and again by Theorem \ref{szego-herm} and \textbf{GLT3} we have
\begin{equation}\label{eq:distrM}
	\left\{\frac{1}{\Delta x}M_{n+1}\right\}_n
	\sim_{{\rm GLT},\sigma,\lambda}
	(\msym,[-\pi,\pi]).
\end{equation}
Therefore, its eigenvalues are $\dfrac{9}{70} d \rho  \left(2 \pm1 \right)  \left( \tfrac{1}{2} \pm \tfrac{1}{16} \right)$.

\begin{figure}[!t]
	\begin{center}
		\begin{tabular}{cc}
			\includegraphics[width=0.4\textwidth]{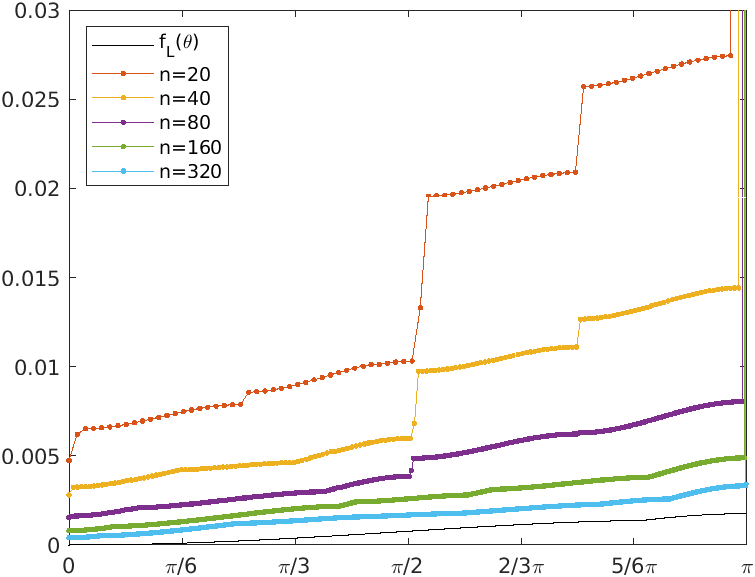} &
			\includegraphics[width=0.4\textwidth]{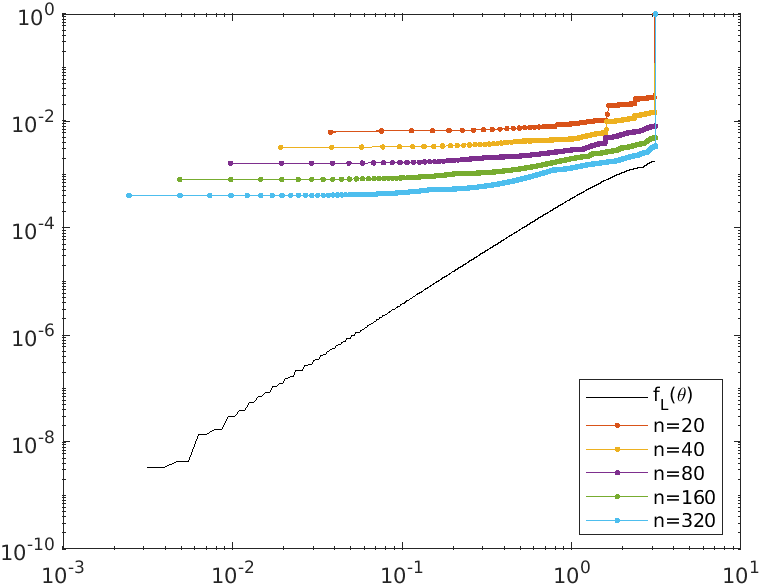}\\
			(a) & (b)
		\end{tabular}
	\end{center}
	\caption{(a)The spectrum of  $(M_{n+1}+L_{n+1})$ with different number of cells vs sampling of the eigenvalue functions of $\lsym(\theta)$ associated to the only matrix $L_{n+1}$; (b) is the same picture, but in bilogarithmic scale.}
	\label{fig:symbolLM}
\end{figure}

In Fig.~\ref{fig:symbolM} we compare an equispaced sampling of the eigenvalues of $\msym(\theta)$ 
with the spectrum of the mass matrix-sequences and we see that the matching is getting better and better as the number of cells increases.

Since the $(1,1)$ block of $\mathcal{A}$ is given by the sum of $L_{n+1}$ and $M_{n+1}$, we are interested in the symbol of $\{N_{n+1}=L_{n+1}+M_{n+1}\}_n$. Let us first note that because of the presence of $\dx$ in its definition, $M_{n+1}$ is a norm-correction of $L_{n+1}$ and that $N_{n+1}$ is real symmetric when boundary conditions are excluded. Then, by using Proposition \ref{prop:zero}, equation \eqref{eq:distrL}, and \textbf{GLT1-4} we have that
\begin{equation}\label{eq:distrN}
	\left\{N_{n+1}\right\}_n\sim_{{\rm GLT},\sigma,\lambda}(\lsym,[-\pi,\pi]).
\end{equation}
Fig. \ref{fig:symbolLM} checks numerically relation \eqref{eq:distrN} by comparing the eigenvalues of $N_{n+1}$ modified by the boundary conditions (see Remark \ref{rem:BCs}) with an equispaced sampling of the eigenvalue functions of $\lsym(\theta)$.

\begin{figure}[!h]
	\begin{center}
	 \includegraphics{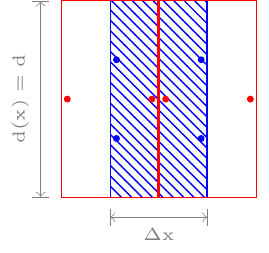}
	\end{center}
	\caption{Illustration of the stencil that refers to the pressure gradient matrix $G_{n+1,n}$.}
	\label{fig:discretG}
\end{figure}

\paragraph{Gradient operator}
The $(1,2)$ block $G$ of $\mathcal{A}$ in \eqref{eq:coeffm} is organized in blocks of rows, each of size $n_u=4$ (blue degrees of freedom in Fig.~\ref{fig:discretG}); in each row there are $2n_p=4$ nonzero elements (red degrees of freedom in Fig.~\ref{fig:discretG}), half of which are associated with the pressure cell intersecting the velocity cell in its left (respectively right) half.
	
\noindent Therefore the gradient matrix is a $4(n+1)\times2n$ rectangular matrix that, excluding boundary conditions, can be written as		
\[
	G_{n+1,n} = \frac{3}{64}d \dt
	\left[
	\begin{array}{ccccccc}\notag
	g_{0}	 & 0      & \cdots & \cdots & \cdots & 0      \\
	g_{1}	 & g_{0}  & 0	     &        & & \vdots \\
	0	     & g_{1}  & g_{0}	 & 0      & & \vdots \\
	\vdots & \ddots & \ddots & \ddots & \ddots & \vdots \\
	\vdots & & 0      & g_{1}  & g_{0}  & 0      \\
	\vdots & & & 0      & g_{1}	& g_{0}  \\
	0\vphantom{\vdots}	     &\cdots& \cdots & \cdots & 0	    & g_{1} 
	\end{array}
	\right]
\]
where 
$g_{0} =
\left[
	\begin{array}{cc}\notag
	3	& 1 \\
	3	& 1 \\
	1	& 3 \\
	1	& 3 \\
	\end{array}
\right]
$ and $g_{1} = -g_{0}$. 

\begin{figure}[!t]
	\begin{center}
		\begin{tabular}{cc}
		\includegraphics[width=0.4\textwidth]{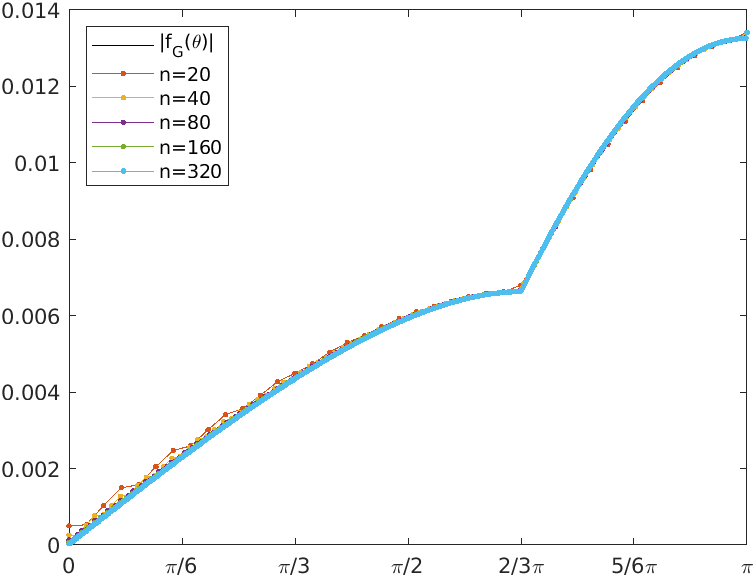} &
		\includegraphics[width=0.4\textwidth]{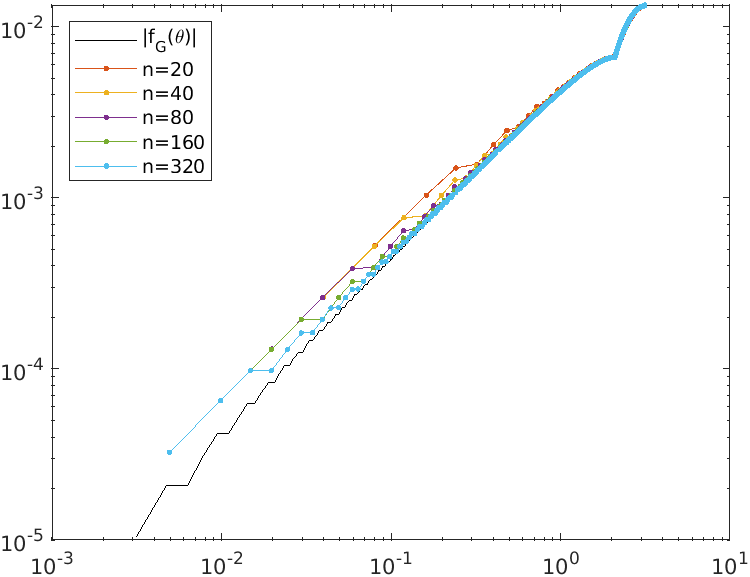}\\
		(a) & (b)
		\end{tabular}
	\end{center}
	\caption{(a) The singular values of $\frac{1}{\dt} G_{n+1,n}$ matrix with different number of cells vs sampling of the singular value functions of $\gsym(\theta)$; (b) is the same picture, but in bilogarithmic scale.}
	\label{fig:symbolG}
\end{figure}

\noindent Similarly to what has been done for the mass matrix-sequence, due to the presence of $\dt$ in $G_{n+1,n}$, we focus on the symbol of the scaled sequence $\{\frac{1}{\dt}G_{n+1,n}\}_n$. Note that $\frac{1}{\dt}G_{n+1,n}$ is a submatrix of a $4\times 2$-block rectangular Toeplitz, precisely $G_{n+1,n}=[T_n(\gsym)]_{n+1,n}$ with 
$\gsym:[-\pi ,\pi] \to \mathbb{C}^{4\times2}$ defined by
\[
	\gsym(\theta)  = \frac{3}{64}d   \left( g_{0} + g_{1} e^{\ii\theta} \right) = \frac{3}{64} d \,  g_{0} (1-e^{\ii\theta}) = - \ii \frac{3}{32} d  \, g_{0}\  e^{\ii\tfrac{\theta}{2}} \sin\left({\tfrac{\theta}{2}}\right),
\]
and thanks to Remark \ref{rem:distr-rect} we deduce

\begin{equation}\label{eq:distrG}
	\left\{\frac{1}{\dt}G_{n+1,n}\right\}_n
	\sim_{\sigma}
	(\gsym,[-\pi,\pi]).
\end{equation}

The singular value decomposition of $g_{0}$ is $U \Sigma V^{T}$ where
\[
	U = \tfrac{1}{2} \left[
	\begin{array}{cccc}\notag
	-1  &  -1 & -1 & -1 \\
	-1  &  -1 &  1 &  1 \\
	-1  &   1 &  1 & -1 \\
	-1  &   1 & -1 &  1 \\
	\end{array} \right]
	\hspace{20pt}
	V = \tfrac{\sqrt{2}}{2} \left[
	\begin{array}{cc}\notag
	-1 &  -1 \\  
	-1  &  1 \\
	\end{array}	\right]
	\hspace{20pt}
	\Sigma = 2\sqrt{2}  \left[
	\begin{array}{cc}\notag
	2  &  0 \\  
	0  &  1 \\
	\end{array}	\right]
\]
and thus the singular value functions of the symbol $\gsym(\theta)$ are $- \frac{3}{8} \sqrt{2}\ii e^{\ii\tfrac{\theta}{2}} \sin\left({\tfrac{\theta}{2}} \right) $ and $- \frac{3}{16} \sqrt{2} \ii e^{\ii\tfrac{\theta}{2}} \sin\left({\tfrac{\theta}{2}}\right)$. Fig.~\ref{fig:symbolG} shows the very good agreement of the spectrum of $\frac{1}{\dt}G_{n+1,n}$ with the sampling of the singular value functions of $\gsym(\theta)$ for different number of cells.

\begin{figure}[!h]
	\begin{center}
			\includegraphics{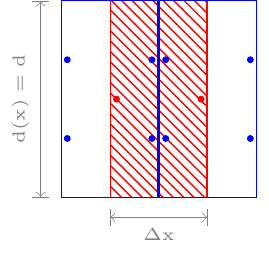}
	\end{center}
	\caption{Illustration of the stencil that refers to the divergence matrix $D_{n,n+1}$.}
	\label{fig:discretD}
\end{figure}

\paragraph{Divergence operator}
The $(2,1)$ block $D$ of the matrix $\mathcal{A}$ is organized in blocks of rows each of size $n_p=2$ (red degrees of freedom in Fig.~\ref{fig:discretD}); in each row there are $2n_u=8$ nonzero elements (blue degrees of freedom in Fig.~\ref{fig:discretD}), half of which are associated with the velocity cell intersecting the pressure cell in its left (respectively right) half.

\begin{figure}[!t]
	\centering
	\begin{tabular}{cc}	
	\includegraphics[width=0.4\textwidth]{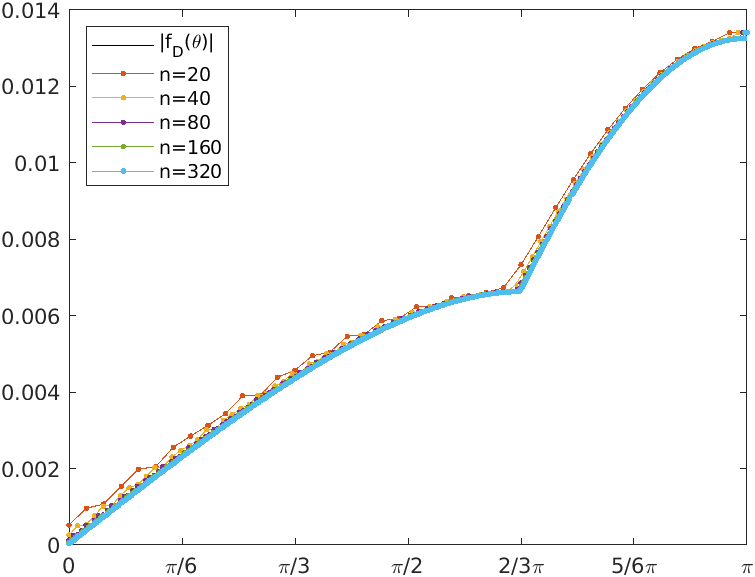} &
	\includegraphics[width=0.4\textwidth]{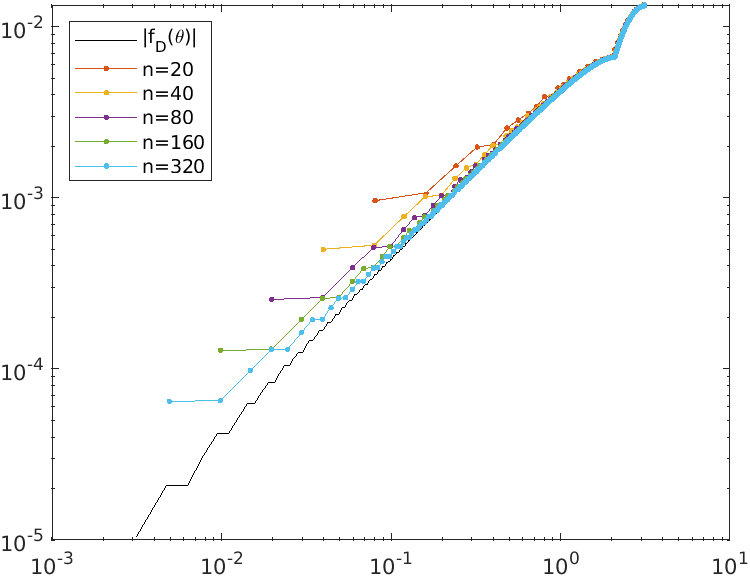}\\
	(a) & (b)
	\end{tabular}
	\caption{(a) The singular values of $D_{n,n+1}$ different number of cells vs sampling of the singular value functions of $\gsym(\theta)$; (b) is the same picture, but in bilogarithmic scale.}
	\label{fig:symbolD}
\end{figure}

Similarly to what we did for the gradient of the pressure, we can define
$d_{0} =
\left[
\begin{array}{cccc}\notag
3 & 3 & 1 & 1 \\
1 & 1 & 3 & 3 \\
\end{array}
\right]=g_0^T
$ and $d_{-1} = -d_{0}$, and we can write the divergence matrix as		
\[
	D_{n,n+1} = \frac{3}{64}d 
	\left[
	\begin{array}{ccccccc}\notag
	d_{0} 	& d_{-1} & 0       & \cdots & \cdots & \cdots &0      \\
	0       & d_{0}  & d_{-1}  & 0      & & & \vdots \\
	\vdots  & 0      & d_{0}   &d_{-1} & 0 &  & \vdots \\
	\vdots  &      &\ddots  & \ddots  & \ddots & \ddots 	     & \vdots \\
	\vdots  &  & & 0       & d_{0}  & d_{-1} & 0      \\
	0\vphantom{\vdots}	      & \cdots & \cdots  & \cdots  &0      & d_{0}	 & d_{-1} \\
	\end{array}
	\right]
\]
Since the matrix $D_{n,n+1}$ is the transpose of $\frac{1}{\Delta t}G_{n+1,n}$, the generating function is 
\[
	\dsym(\theta) = (\gsym(\theta))^*= \ii\frac{3}{32} \, d \, g^T_{0} \ e^{-\ii\tfrac{\theta}{2}}\ \sin\left({\tfrac{\theta}{2}}\right) 
\]
which admits the same singular value functions of $\gsym(\theta)$. Therefore, by Remark \ref{rem:distr-rect} we find 

\begin{equation}\label{eq:distrD}
	\{D_{n,n+1}\}_n\sim_{\sigma}(\dsym,[-\pi,\pi]).
\end{equation}
A comparison of the sampling of the singular values of $\dsym(\theta)$ with the singular values of $D_{n,n+1}$ is shown in Fig.~\ref{fig:symbolD}.

\begin{figure}[!t]
	\centering
	\begin{tabular}{cc}	
	\includegraphics[width=0.4\textwidth]{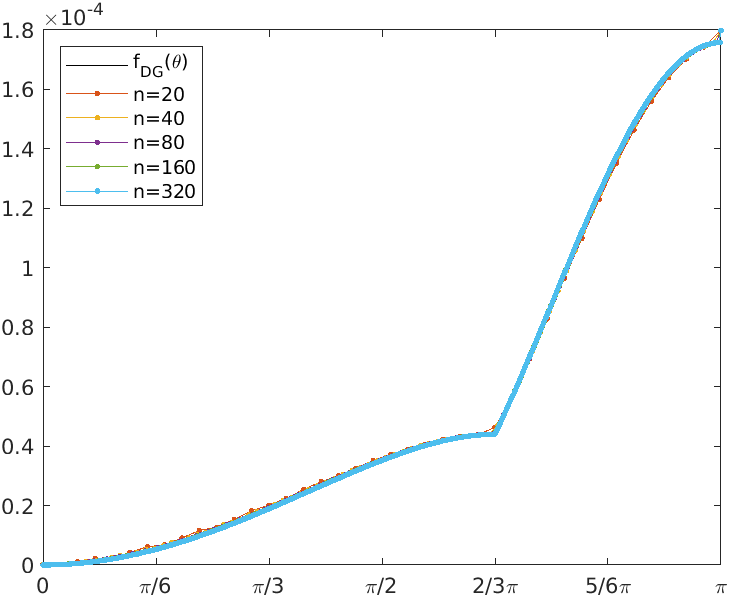} &
	\includegraphics[width=0.4\textwidth]{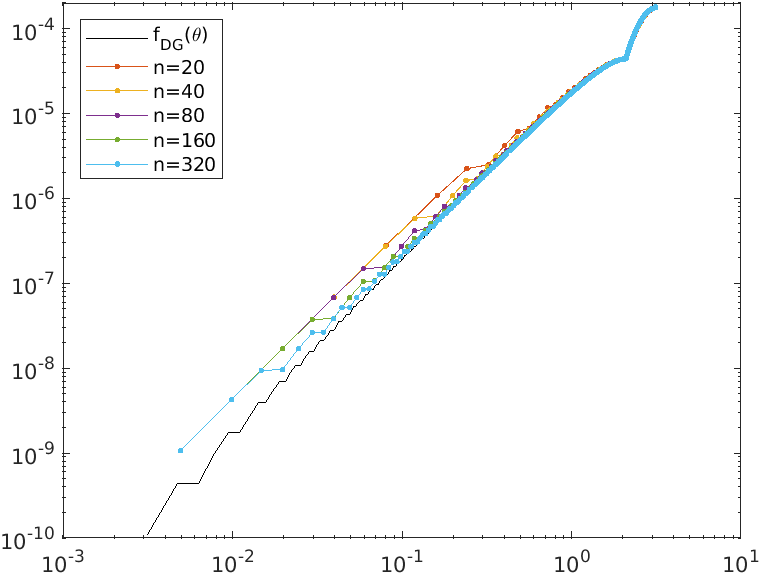}\\
	(a) & (b)
	\end{tabular}
	\caption{(a) The spectrum of the matrix product $\frac{1}{\dt}D_{n,n+1}G_{n+1,n}$ with different number of cells vs sampling of the eigenvalues of $\dsym(\theta)\gsym(\theta)$; (b) is the same picture, but in bilogarithmic scale.}
	\label{fig:symbolDG}
\end{figure} 

\begin{remark}\label{rem:sym_dg}
If we analyse the product of the symbols for $D_{n,n+1}$ and $\frac{1}{\Delta t}G_{n+1,n}$, we obtain a $\C^{2\times2}$-valued symbol:
\begin{align*}
\dsym(\theta)\gsym(\theta)
	= V\Sigma U^{T} U \Sigma V^{T}&= 
	\left[
	\begin{array}{cc}\notag
	5 & 3 \\
	3 & 5 \\
	\end{array}
	\right]  4 \sin^2\left(\tfrac{\theta}{2}\right)
	\left(\frac{3}{32}d \right)^2\\
	&=\left[
	\begin{array}{cc}\notag
	5 & 3 \\
	3 & 5 \\
	\end{array}
	\right]  2 (1-\cos\theta)
	\left(\frac{3}{32}d \right)^2 
\end{align*}
Its eigenvalue functions are $4(1-\cos\theta)\left(\frac{3}{64}d \right)^2$ and  $16 (1-\cos\theta)\left(\frac{3}{64}d \right)^2$. Notice that, since $D_{n,n+1}=[T_n(\dsym)]_{n,n+1}$ and $\frac{1}{\Delta t}G_{n+1,n}=[T_n(\gsym)]_{n+1,n}$, then $\frac{1}{\Delta t}D_{n,n+1}G_{n+1,n}$ is a principal submatrix of $T_{n}(\dsym)T_{n}(\gsym)$. Therefore, thanks to Theorem \ref{th:distr-rect} and Remark \ref{rem:th-ext}, $\dsym(\theta)\gsym(\theta)$ is the spectral symbol of $\{T_n(\dsym)T_n(\gsym)\}_n$ and, by Theorem \ref{th:extradimensional}, it is also the symbol of $\{\frac{1}{\Delta t}D_{n,n+1}G_{n+1,n}\}_n$. As a consequence, we expect that a sampling of the eigenvalue functions of $\dsym(\theta)\gsym(\theta)$ provides an approximation of the spectrum of $\frac{1}{\Delta t}D_{n,n+1}G_{n+1,n}$. This is confirmed by  
Fig.~\ref{fig:symbolDG}.
\end{remark}

\paragraph{Penalty term for pressure}

\begin{figure}[!t]
	\begin{center}
		\includegraphics[width=0.4\textwidth]{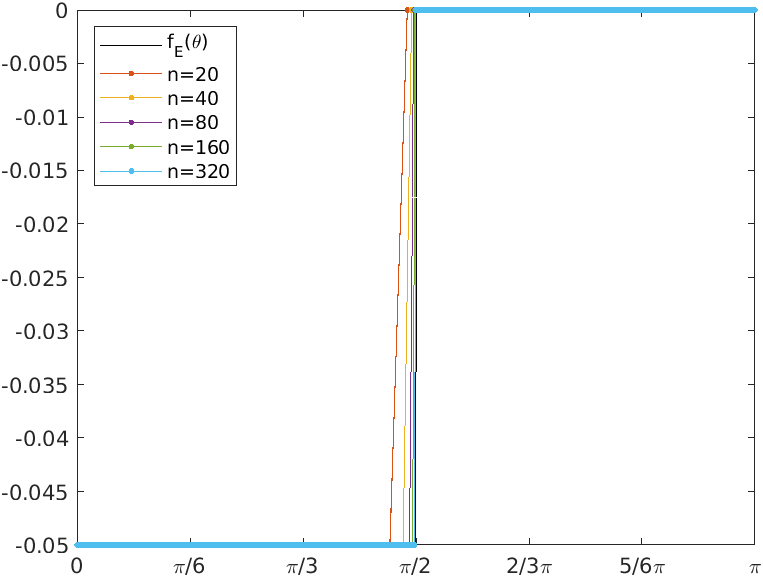}
	\end{center}
	\caption{The spectrum of $ \frac{1}{\dx} E_n$ with different number of cells vs sampling of the eigenvalue functions of $\esym(\theta)$.} 
\end{figure}

The $(2,2)$ block of  matrix $\mathcal{A}$ is organized in blocks of rows, each of size $n_p=2$ and it has the following form
\[
E_n = d \dx
\tridiag
\left[
\begin{array}{cc|cc|cc}\notag
0 & 1 & -1 & 0 &  0 & 0 \\
0 &  0 & 0 & -1 & 1 & 0 \\
\end{array}
\right],
\]
where $n$ is the number of pressure cells.
The symbol associated to the scaled matrix-sequence $\{\frac{1}{\Delta x}E_n\}_n$ is the function $\esym:[-\pi ,\pi] \to \mathbb{C}^{2\times2}$  and can be written as
\[
\esym(\theta)  = 
d
\begin{bmatrix}
-1 & e^{\ii\theta} \\
e^{-\ii\theta} & -1 \\
\end{bmatrix}
\]
and so its eigenvalues are $0$ and $-2d$, while its eigenvectors are
$
\begin{pmatrix}
e^{\ii\theta}\\
\ii
\end{pmatrix}$ 
and
$\begin{pmatrix}
-e^{\ii\theta}\\
\ii
\end{pmatrix}.$ 
Since $E_n$ is real symmetric, by \textbf{GLT3} and \textbf{GLT1} we obtain 
\begin{equation}\label{eq:distrE}
\left\{\frac{1}{\dx}E_n\right\}_n\sim_{{\rm GLT},\sigma,\lambda}(\esym,[-\pi,\pi]).
\end{equation}

\subsection{Spectral study of the Schur complement}\label{sub:schur}

\begin{figure}[!t]
	\begin{center}
		\begin{tabular}{cc}	
			\includegraphics[width=0.4\textwidth]{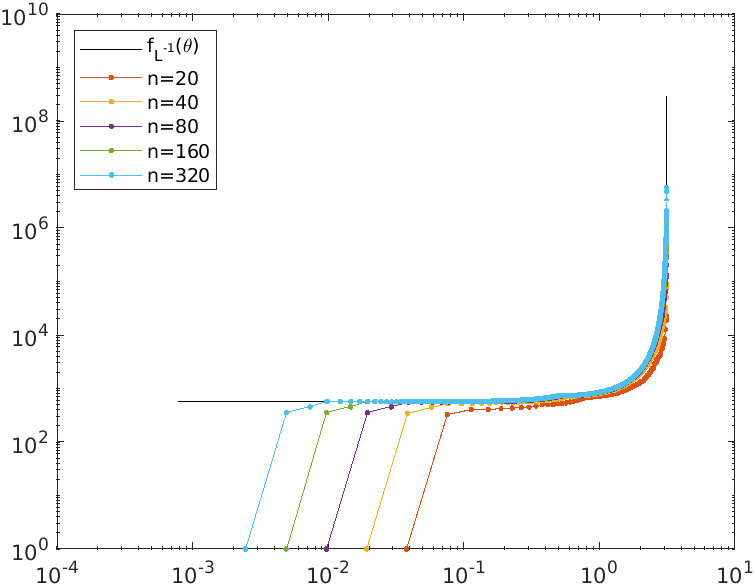} &
			\includegraphics[width=0.4\textwidth]{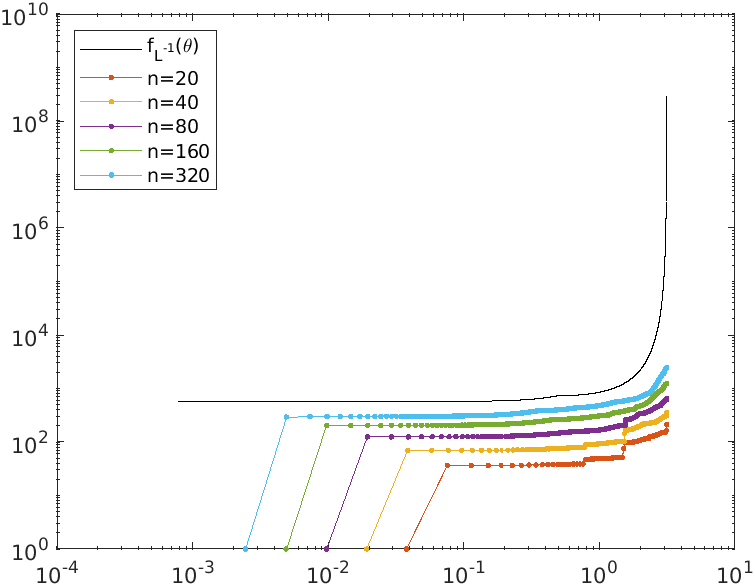}\\
			(a) $L_{n+1}^{-1}$ & (b) $(L_{n+1}+M_{n+1})^{-1}$
		\end{tabular}
	\end{center}
	\caption{The spectrum of $L_{n+1}^{-1}$ and $(L_{n+1}+M_{n+1})^{-1}$
		vs the eigenvalue functions of $\lsym^{-1}(\theta)$.} 
	\label{fig:symbolLInv}
\end{figure}

We now study the spectral distribution of the Schur complement of $\mathcal{A}$. The formal expression of the Schur complement involves inversion of the $(1,1)$ block of the matrix system and the multiplication by the $(1,2)$ and $(2,1)$ blocks that is: $S_n = E_n -D_{n,n+1}N_{n+1}^{-1}G_{n+1,n}$. To compute the symbol of the Schur complement sequence  we need to compute the symbol of $\{(L_{n+1}+M_{n+1})^{-1}\}_n$. Thanks to relation \eqref{eq:distrN} and to \textbf{GLT1-2} we have 		\begin{equation}\label{eq:distrLinv}
\{(L_{n+1}+M_{n+1})^{-1}\}_n\sim_{\lambda}(\lsym^{-1},[-\pi,\pi])
\end{equation}
with
\[
\lsym^{-1}(\theta)  =  \frac{b}{1-cos\theta}
\begin{bmatrix}
8 & 1 & 0 & 0 \\
1 & 8 & 0 & 0 \\
0 & 0 & 8 & 1 \\
0 & 0 & 1 & 8 \\
\end{bmatrix}	
\]
where $b = \frac{560}{1701} \frac{1}{\mu d c}$. $\lsym^{-1}$ has two eigenvalue functions $\frac{9b}{1-cos\theta}$ and $\frac{7b}{1-cos\theta}$, each with multiplicity 2. 
Following \eqref{eq:distrLinv}, in Fig.~\ref{fig:symbolLInv} we compare the spectrum of $L_{n+1}^{-1}$ and of $(L_{n+1}+M_{n+1})^{-1}$ with a sampling of the eigenvalue functions of $\lsym^{-1}(\theta)$. In both cases the spectrum of the matrix has the same behavior of the symbol.

\begin{figure}[!t]
	\begin{center}
		\begin{tabular}{cc}	
		\includegraphics[width=0.4\textwidth]{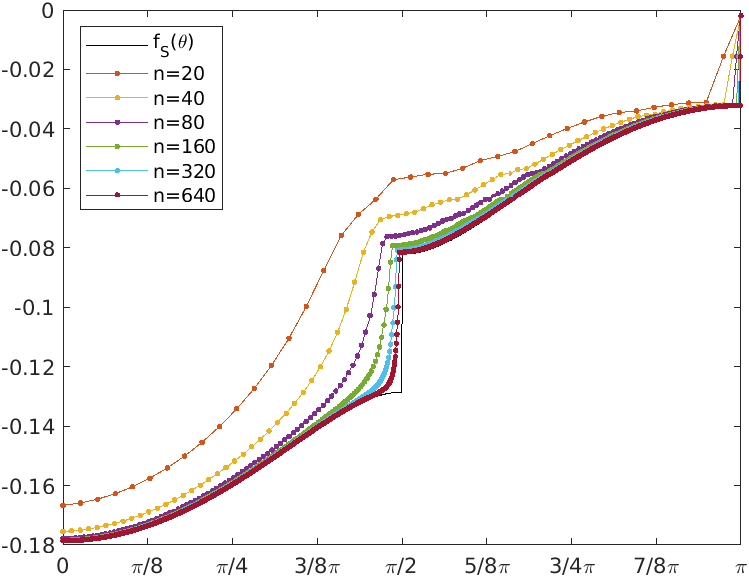} &
		\includegraphics[width=0.4\textwidth]{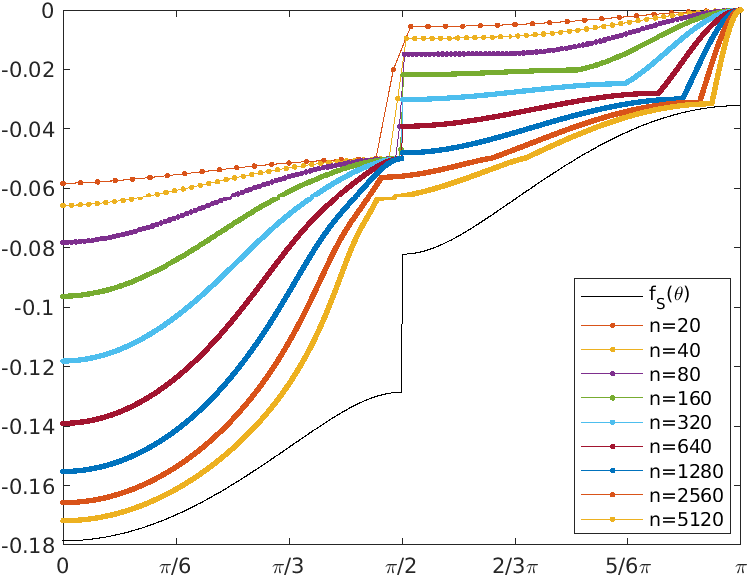}\\
		(a) $E_{n}- D_{n, n+1} L_{n+1}^{-1} G_{n+1, n}$ & (b) $E_{n}- D_{n, n+1} (L+M)_{n+1}^{-1} G_{n+1, n}$
		\end{tabular}
	\end{center}
	\caption{The spectrum of the matrix $\frac{1}{\dt}S_n $ with different number of cells vs sampling of the eigenvalue functions of the symbol $\ssym(\theta)$ In (a), the (1,1) block contains only the $L_{n+1}$ term, while in (b) the block $N_{n+1}$ contains $L_{n+1}+M_{n+1}$.}
	\label{fig:symbolS}

\end{figure}

\begin{figure}[!t]
	\begin{center}
		\includegraphics[width=.7\textwidth]{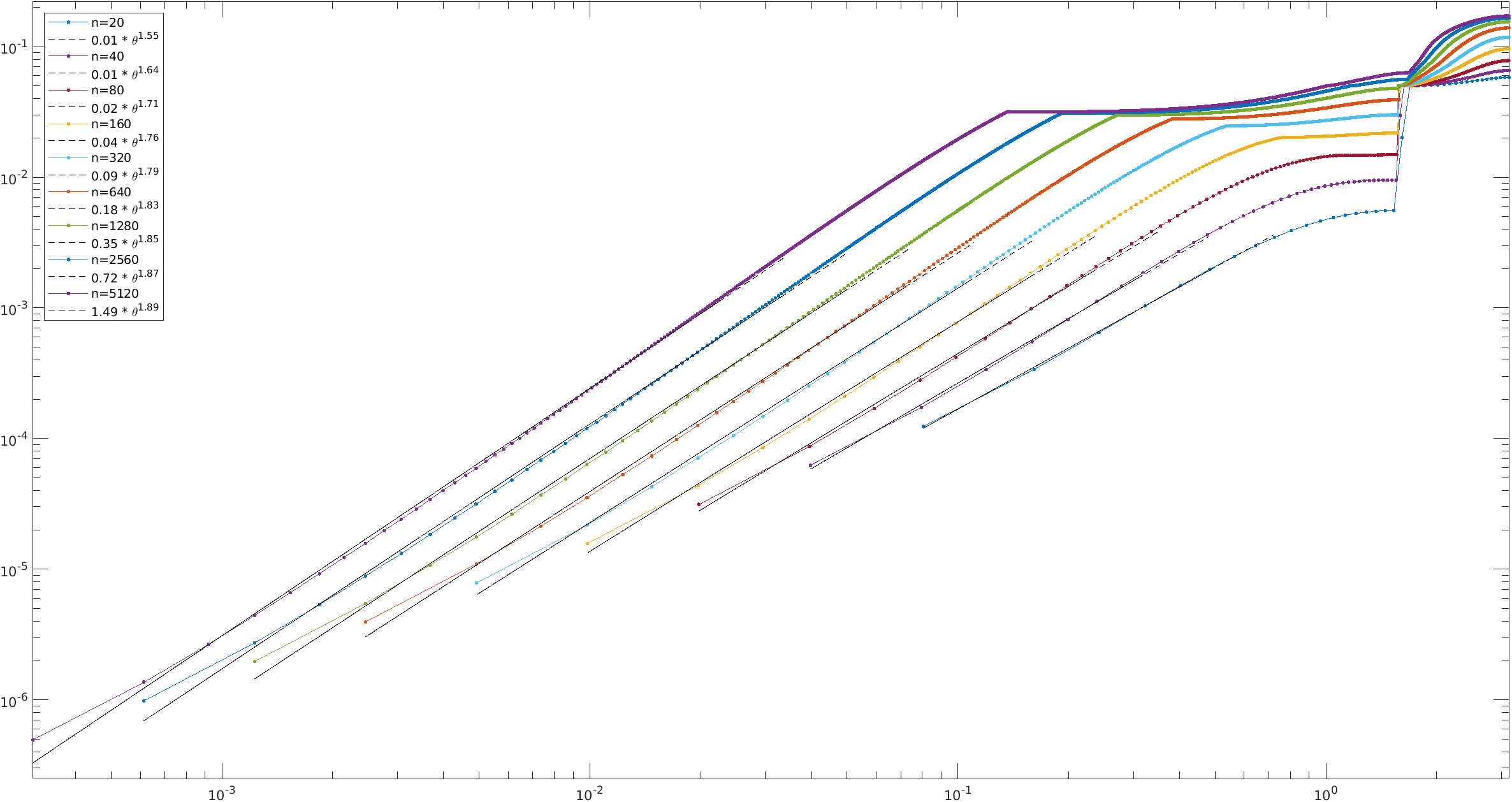}
	\end{center}
	\caption{Smallest eigenvalues of $-\frac{1}{\dt}S_n$ and best fits with functions of the type $c\cdot\theta^{\gamma}$: for large $n$ the order $\gamma$ is, as expected, approximately 2.}
	\label{fig:symbolSwithM}
\end{figure}

At this point we can focus on the symbol of a properly scaled Schur complement sequence: $\{\frac{1}{\dt}S_n\}_n$. We know that $\frac{1}{\dt}S_n$ is a principal submatrix of 
$$
\tilde{S}_{n}:=T_{n}\left(\frac{1}{c}\esym\right)-T_n(\dsym)T_n(\lsym)^{-1}T_n(\gsym)+Z_n
,
$$
$Z_n$ being a correction-term. Since we are assuming that $c=\frac{\dt}{\dx}=\Ogrande(1)$ and since $\lsym(\theta)$ is an Hermitian positive definite matrix-valued function, by combining Theorem \ref{th:distr-rect-bis}, and equations \eqref{eq:distrG}, \eqref{eq:distrD}, \eqref{eq:distrE}, \eqref{eq:distrLinv} it holds that 
\begin{equation*}
\left\{T_{n}\left(\frac{1}{c}\esym\right)-T_n(\dsym)T_n(\lsym)^{-1}T_n(\gsym)\right\}_n\sim_{\sigma,\lambda}(\ssym,[-\pi,\pi])
\end{equation*}
where
\[
\ssym(\theta)  = \frac{1}{c}\esym(\theta)-\dsym(\theta)\lsym^{-1}(\theta)\gsym(\theta)= \frac{d}{c}
\begin{bmatrix} 
-1 - 5 \tfrac{a}{\mu}            & e^{\ii\theta} - 3 \tfrac{a}{\mu}  \\
e^{-\ii\theta} - 3\tfrac{a}{\mu} & -1 - 5 \tfrac{a}{\mu} \\
\end{bmatrix}	
\]
and $a = \tfrac{105}{2016}$. This combined with Theorem \ref{thm:quasi-herm}  guarantees that 
\begin{equation*}
\left\{\tilde{S}_{n}\right\}_n\sim_{\lambda}(\ssym,[-\pi,\pi])
\end{equation*}
and consequently
\begin{equation}\label{eq:distrS}
\left\{\frac{1}{\dt}S_n\right\}_n\sim_{\lambda}(\ssym,[-\pi,\pi]).
\end{equation}
The eigenvalue functions of $\ssym(\theta)$ are $ \frac{d}{c} \left(-1 - 5 \tfrac{a}{\mu} \pm \sqrt{1 + 9 \tfrac{a^2}{\mu^2}  -6\tfrac{a}{\mu}cos\theta}\right)$. 
In Fig.~\ref{fig:symbolS} we compare a sampling of the eigenvalue functions of $\ssym(\theta)$ with the spectrum of $\frac{1}{\dt}S_n$ for different grid refinements. In the right panel, we consider the complete matrix $\mathcal{A}$ with $N_{n+1}=L_{n+1}+M_{n+1}$, while in the left panel we show the situation when replacing $N_{n+1}$ with $L_{n+1}$. 
Moreover, in Fig.~\ref{fig:symbolSwithM} we compare the minimal eigenvalues of 
$-\frac{1}{\dt}S_n$ with functions of type $c\cdot\theta^{\gamma}$ and we see that for large $n$ the order $\gamma$ is approximately 2.

\begin{remark}
We stress that, thanks to the newly introduced Theorem \ref{th:distr-rect-bis}, computing the symbol of the product $D_{n,n+1}N_{n+1}^{-1}G_{n+1,n}$ immediately follows by using standard spectral distribution tools as Theorem \ref{thm:quasi-herm}. The same result could be obtained following the much more involved approach used in \cite{system-PDE-FEM}. Such approach asks to first extend the rectangular matrices $D_{n,n+1}$, $G_{n+1,n}$ to proper square block Toeplitz matrices, and then use the GLT machinery to compute the symbol of their product with $N_{n+1}^{-1}$. Finally, the symbol of the original product is recovered by projecting on the obtained matrix through ad hoc downsampling matrices and by leveraging the results on the symbol of projected Toeplitz matrices designed in the context of multigrid methods \cite{serra2004multigrid}.
\end{remark}

\begin{figure}[!h]
	\begin{center}
		\begin{tabular}{cc}
		\includegraphics[width=0.40\textwidth]{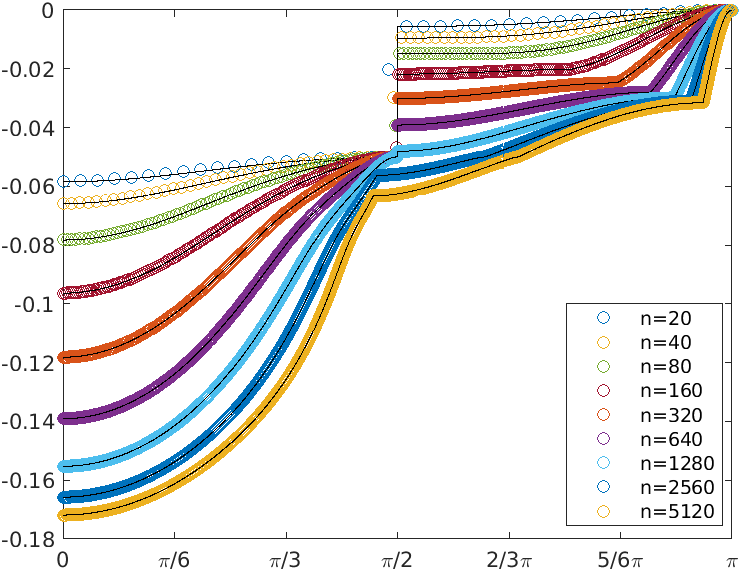}&
		\includegraphics[width=0.40\textwidth]{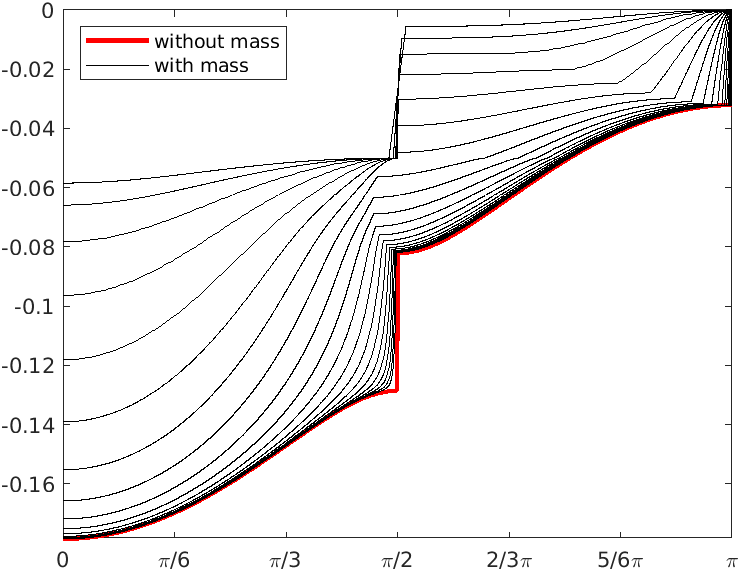}	\\
		(a) & (b)
		\end{tabular}
 
	\end{center}
	\caption{(a) The spectrum of the matrix $\frac{1}{\dt }S_n$ with different number of cells vs sampling of the eigenvalues of $\ssym_{\dx}(\theta)$, (b) Visual convergence of the generating function $\ssym_{\dx}(\theta)$ (black lines) to $\ssym(\theta)$ (red line) as $\dx\rightarrow 0$.}
	\label{fig:symbolSwithM:fSM}
\end{figure}

Aside from the symbol $\ssym(\theta)$, having in mind to build a preconditioner for the Schur matrix, we compute also the generating function of $\frac{1}{\dt}S_n$ for a fixed $n$, that is for a fixed $\dx$. Here we keep the contribution of the mass matrix in $N_{n+1}$. As a result, we get
\begin{equation}\label{eq:sym_deltax}
\ssym_{\dx}(\theta)  = \tfrac{d}{c}
\begin{bmatrix}
-1 - (5a(\theta)-3\dx\rho)b(\theta)c             & e^{\ii\theta} - (3a(\theta)-5\dx\rho)b(\theta)c  \\
e^{-\ii\theta} - (3a(\theta)-5\dx\rho)b(\theta)c & -1 - (5a(\theta)-3\dx\rho)b(\theta)c \\
\end{bmatrix}	
\end{equation}
with  $a(\theta) = 6\left(1-cos\theta \right) \mu c +2\dx\rho$ and $b(\theta) = \frac{315}{1008}\frac{\left(1-cos\theta \right)}{a(\theta) ^2 - \dx^2 \rho^2}$. As shown in Fig.~\ref{fig:symbolSwithM:fSM}(a), the sampling of the eigenvalue functions of $\ssym_{\dx}(\theta)$ perfectly matches the spectrum of the corresponding Schur matrix, and this paves the way to design a preconditioner that instead of $\ssym(\theta)$ involves $\ssym_{\dx}(\theta)$. Of course, in the limit when $\dx$ goes to zero, the symbol is equal to $\ssym(\theta)$. As a confirmation see Fig.~\ref{fig:symbolSwithM:fSM}(b).

\subsection{Spectral study of the coefficient matrix}
\begin{figure}[!t]
\begin{center}
	\begin{tabular}{cc}
		\includegraphics[width=0.4\textwidth]{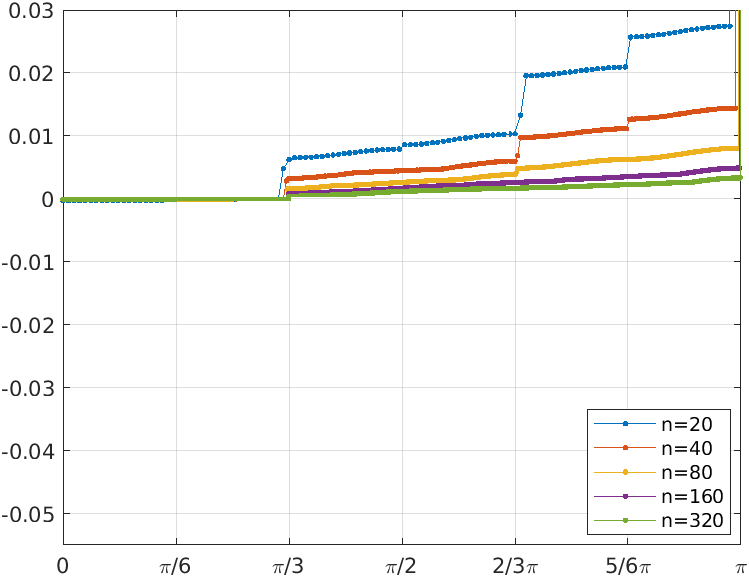} &
		\includegraphics[width=0.4\textwidth]{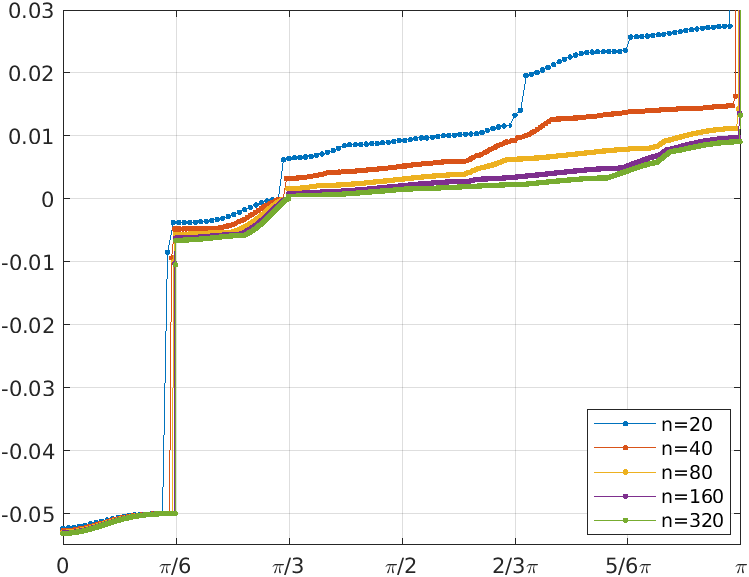}\\
		(a) original $\mathcal{A}$ & (b) scaled $\mathcal{A}V$
	\end{tabular}
\end{center}
\caption{The spectrum of the coefficient matrix.}
\label{fig:scaleA}
\end{figure}

The results obtained in Subsections \ref{sub:blocks}-\ref{sub:schur} suggest to scale the coefficient matrix $\mathcal{A}$ by columns through the following matrix
\[ V=\begin{bmatrix}\notag
I  & 0 \\
0  & \frac{1}{\dt} I \\
\end{bmatrix},
\]
that is to solve the system $\mathcal{A}_n\textbf{x} = \textbf{f}$, with $\mathcal{A}_n:=\mathcal{A}V$  in place of system \eqref{eq:coeffm}. As a result of the scaling, the blocks $\frac{1}{\dt}G_{n+1,n}$ and $\frac{1}{\dt}E_n$ of $\mathcal{A}_n$ have size $O(1)$, similar to the size of $N_{n+1}$ and $D_{n,n+1}$, which remain unchanged. Moreover, the scaling improves the arrangement of the eigenvalues of $\mathcal{A}$ since the small negative eigenvalues are shifted towards negative values of larger modulus, as we can see in Fig.~\ref{fig:scaleA}. Indeed, excluding the boundary conditions and due to the block-factorization  \[ \mathcal{A}_n=WDW^T=\begin{bmatrix}\notag
I_{n+1}  & 0 \\
D_{n,n+1}N_{n+1}^{-1}  & I_n\\
\end{bmatrix}
\begin{bmatrix}\notag
N_{n+1}  & 0 \\
0  & \frac{1}{\dt}S_n  \\
\end{bmatrix}
\begin{bmatrix}\notag
I_{n+1}  &N_{n+1}^{-1}\frac{1}{\dt} G_{n+1,n}  \\
0  & I_n \\
\end{bmatrix},
\]
by the Sylvester inertia law we can infer that the signature of $\mathcal{A}_n$ is the same of the signature of the diagonal matrix formed by $N_{n+1}$ and $\frac{1}{\dt}S_n=\frac{1}{\dt}(E_n-D_{n,n+1}N^{-1}_{n+1}G_{n+1,n})$, which we know has negative eigenvalues distributed according to $\ssym(\theta)$.

\begin{figure}[!t]
	\begin{center}
		\begin{tabular}{cc}
			\includegraphics[width=0.4\textwidth]{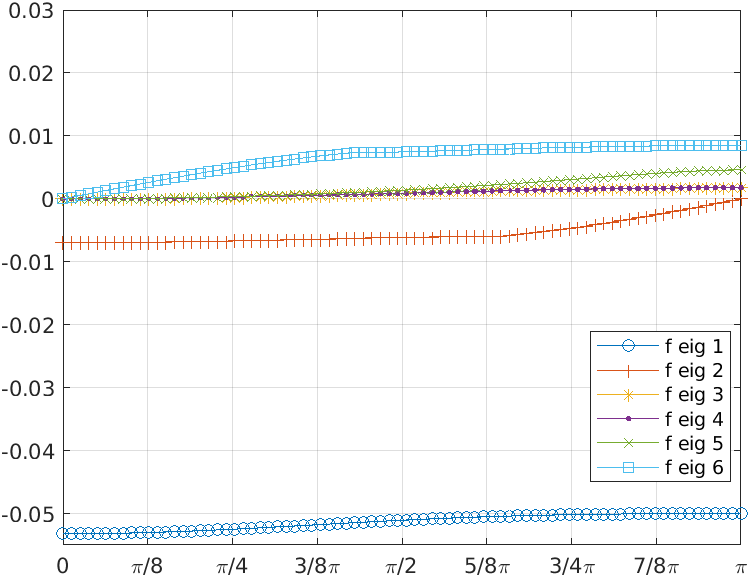} &
			\includegraphics[width=0.4\textwidth]{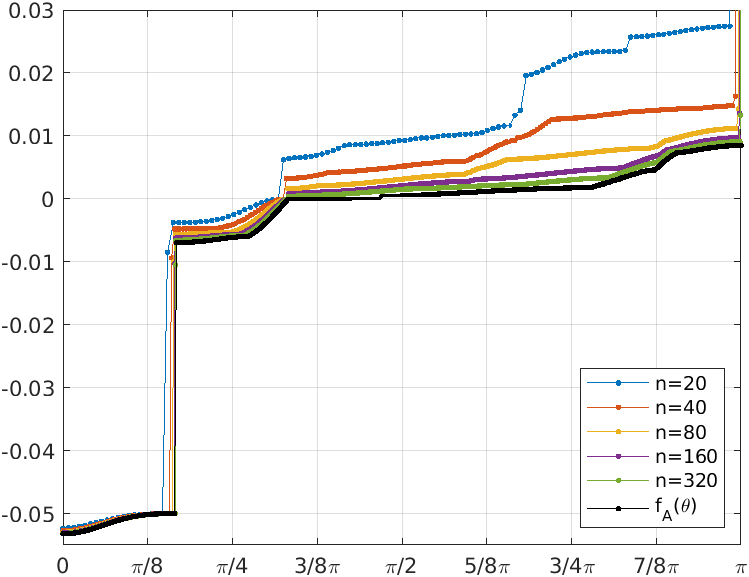}\\
			(a)  & (b) 
		\end{tabular}
	\end{center}
	\caption{(a) A plot of the eigenvalue functions of $\mathscr{F}(\theta)$ made without knowing their analytical expression, (b) The spectrum of the scaled coefficient matrix $\mathcal{A}V$ with different number of cells vs the sampling of the eigenvalue functions of $\mathscr{F}(\theta)$.}
	\label{fig:scaleA_comp}
\end{figure}

In order to obtain the symbol of $\{\mathcal{A}_n\}_n$, let us observe that, when including also the boundary conditions, $\mathcal{A}_n= \tilde{\mathcal{A}}_n+\mathcal{Q}_n$, where $\tilde{\mathcal{A}}_n$ is Hermitian and  $\mathcal{Q}_n$ is a correction term. Let us observe that $\tilde{\mathcal{A}}_n$ is a principal submatrix (obtained removing the last $2$ rows and the last $2$ columns) of the matrix
\begin{align*}
\mathcal{B}_n:&={\begin{bmatrix}
T_{n}(\lsym)+\dx T_n(\msym) & T_{n}(\gsym)\\
T_{n}(\dsym) & T_{n}(\frac{1}{c}\esym)
\end{bmatrix}}\\
&=\begin{bmatrix}
T_{n}(\lsym) & T_{n}(\gsym)\\
T_{n}(\dsym) & T_{n}(\frac{1}{c}\esym)
\end{bmatrix}+\dx\begin{bmatrix}
T_n(\msym) & O\\
O& O
\end{bmatrix}.
\end{align*}
Now, by Theorem \ref{th:wholeA}, the two involved matrices are similar that is
$$\mathcal{B}_n\sim T_n(\mathscr{F})+\dx T_n(\mathscr{C})$$
with $\mathscr{F}:=\begin{bmatrix}
\lsym & \gsym\\
\dsym & \frac{1}{c}\esym
\end{bmatrix}$ and $\mathscr{C}:=\begin{bmatrix}
\msym & 0\\
0 & 0
\end{bmatrix}$.
Therefore, 
\begin{equation*}
\{\mathcal{B}_n\}_n\sim_{\lambda}(\mathscr{F},[-\pi,\pi]),
\end{equation*}
and this, thanks to Theorem \ref{th:extradimensional}, implies that
\begin{equation*}
\{\tilde{\mathcal{A}}_n\}_n\sim_{\lambda}(\mathscr{F},[-\pi,\pi]).
\end{equation*}
Finally, by following the same argument applied in the computation of the Schur complement symbol at the beginning of Section \ref{sub:schur}, by using again Theorem \ref{thm:quasi-herm} we arrive at

\begin{equation*}
\{\mathcal{A}_n\}_n\sim_{\lambda}(\mathscr{F},[-\pi,\pi]).
\end{equation*}
Since the symbol $\mathscr{F}$ is a $6\times6$ matrix-valued function, retrieving an analytical expression for its eigenvalue functions asks for some extra computation, but we can easily give a numerical representation of them which is sufficient for our aims simply following these three steps:

\begin{itemize}
	\item evaluate the symbol $\mathscr{F}$ on an equispaced grid in $[0,\pi]$;
	\item for each obtained $6\times6$ matrix compute the spectrum;
	\item take all the smallest eigenvalues as a representation of $\lambda_1(\mathscr{F})$ and so on so forth till the largest eigenvalues as a representation of $\lambda_6(\mathscr{F})$.
\end{itemize}
Fig. \ref{fig:scaleA_comp}(a) has been realized following the previous steps. Notice that two eigenvalue functions of $\mathscr{F}$ show the same behavior and we suspect they indeed have the same analytical expression. Fig. \ref{fig:scaleA_comp}(b) compares the equispaced sampling of the eigenvalue functions with the actual eigenvalues of the coefficient matrix and highlights an improving matching as the matrix-size increases.

\begin{remark}\label{rem:general}
The eigenvalue structure in the general case of a variable cross-section $d=d(x)$ does not pose technical problems and in reality it is perfectly covered by the GLT theory: more specifically, we refer to item \textbf{GLT1} where the GLT symbol depends on $(x,\theta)\in [0,1]\times [-\pi,\pi]$ and where $x$ is in our context exactly the scaled physical variable of the coefficient $d=d(x)$.

The case of a variation of the degrees $n_x$, $n_y$ is more delicate to treat, since, in this setting, the size of the basic small blocks of the matrix is affected. This is the parameter $s$ defining the range $\mathbb{C}^{s\times s}$ of the symbol $\kappa$ in the GLT theory (see Section \ref{sec:prelim}). Despite the theoretical difficulty of treating a varying parameter $s$ for a precise spectral analysis, as shown in the next section, the performances of our preconditioning techniques are satisfactory also in this tricky setting.
\end{remark}

\begin{remark}\label{rem:3D}
Our discretization can be extended to three-dimensional pipes by introducing tensor product shape functions in the transverse plane, using polynomial degrees $n_{y}$ and $n_{z}$ for the velocity. 
Leaving fixed $n_x=1$ for the pressure variable, our theory should extend to this more general setting and yield a symbol for the $(1,1)$-block of the coefficient matrix with  values in $\C^{2(n_y-1)(n_z-1)\times2(n_y-1)(n_z-1)}$,
symbols for $(1,2)$- and $(2,1)$-blocks in
$\C^{2(n_y-1)(n_z-1)\times2}$ and $\C^{2\times2(n_y-1)(n_z-1)}$ respectively.
In any case, the symbol for $(2,2)$-block and the Schur complement will still take values in $\C^{2\times2}$ independently of $n_y$ and $n_z$. The size $2\times2$ for the symbol of the Schur complement is controlled by the choice of $n_x=1$ for the pressure variable, 
and for larger $n_x$ the symbol of the Schur complement should take values in
$\C^{(n_x+1)\times(n_x+1)}$.
\end{remark}


\section{Numerical experiments}
\label{sec:numer}

In this section we focus on the solution of system \eqref{eq:coeffm} by leveraging the spectral findings in \S\ref{sec:spectr} and with the help of the PETSc \cite{petsc-efficient,petsc-user-ref} library. To ease the notation, here after we omit the subscripts for the blocks $N_{n+1},G_{n+1,n},D_{n,n+1},E_{n}$ of $\mathcal{A}$. The main solver for $\mathcal{A}_n=\mathcal{A}V$, say $\mathcal{K}_{\mathcal{A}}$, is GMRES and the preconditioner of this Krylov solver is based on the Schur complement; more precisely, an application of the preconditioner consists in solving 
\[
	\hat{S} \hat{p} =  r_p - D \widetilde{N^{-1}} r_u
	\qquad
	\hat{u} = \widetilde{N^{-1}} (r_u-\frac{1}{\dt}Gr_p)
\]
where the block vector $\left(\begin{smallmatrix}r_u\\r_p\end{smallmatrix}\right)$
is the residual. 

If the inversion of $N$ was exact and $\hat{S}$
was the exact Schur complement of $\mathcal{A}_n$, the main solver $\mathcal{K}_{\mathcal{A}}$ would of course be a direct method.
Here above, instead, $\widetilde{N^{-1}}$ denotes the application of a suitable Krylov solver, say $\mathcal{K}_{N}$, to the linear operator $N$ and in our numerical experiments this was chosen as GMRES with a relative stopping tolerance $10^{-5}$ and ILU(0) preconditioner, since $N$ is a narrow-banded matrix.
Further, the Schur complement is approximated by $\hat{S}=\tfrac{1}{\dt}(E-D\widetilde{N^{-1}}G)$. However, since the inverse of $N$ is approximated by the action of the solver $\mathcal{K}_N$, matrix $\hat{S}$ cannot be explicitly assembled, although its action on any vector can be computed with a call to $\mathcal{K}_N$. 

The solution of the system with matrix $\hat{S}$ required in the preconditioner inside $\mathcal{K}_{\mathcal{A}}$ is then performed with a Krylov solver, say $\mathcal{K}_{\hat{S}}$.  In $\mathcal{K}_{\hat{S}}$, the matrix-vector multiplication is performed as described above, while the preconditioner is the block circulant preconditioner
generated by $\ssym_{\dx}(\theta)$ given in \eqref{eq:sym_deltax}, that is (see Theorem \ref{circ-s})
\begin{equation*}
C_n(\ssym_{\dx})=(F_n\otimes I_2) D_n(\ssym_{\dx}) (F_n^*\otimes I_2)
\end{equation*}
with 
\begin{equation*}
D_n(\ssym_{\dx})={\rm diag}_{r=0,\ldots,n-1}(\ssym_{\dx}(\theta_r)), \quad F_n=\frac{1}{\sqrt{n}}\left[e^{-\ii j\theta_r}\right]_{j,r=0}^{n-1}, \quad \theta_r=\frac{2\pi r}{n}.
\end{equation*}
More precisely, since $\ssym_{\dx}(\theta)$ has a unique zero eigenvalue at $\theta_0=0$, we use as preconditioner
\begin{equation}\label{eq:circ_prec}
\mathcal{C}_n:=C_n(\ssym_{\dx})+\frac{1}{(2n)^2}{\bf 1}^T{\bf 1}\otimes\begin{bmatrix}1&1\\1&1\end{bmatrix}
\end{equation}
with ${\bf 1}=[1,\ldots,1]\in\mathbb{R}^{n}$, that is we introduce a circulant rank-one correction aimed at avoiding singular matrices.
We notice that $\{\mathcal{C}_n\}_n$ and the sequence of the Schur complements are GLT matrix-sequences having the same symbol, i.e., $\ssym(\theta)$. Therefore, since $\ssym(\theta)$ is not singular by \textbf{GLT2} we infer that the sequence of the preconditioned matrices is a GLT with symbol 1. Given the one-level structure of the involved matrices, we expect that the related preconditioned Krylov solvers converge within a constant number of iterations independent of the matrix-size, just because the number of possible outliers is bounded from above by a constant independent of the mesh-size. Hence the global cost is given by $O(n\log n)$ arithmetic operations when using the standard FFT based approach for treating the proposed block circulant preconditioner. Furthermore it is worth mentioning that reduction to the optimal cost of $O(n)$ arithmetic operations is possible by using specialized multigrid solvers designed ad hoc for circulant structures \cite{serra2004multigrid}.

The circulant preconditioner is applied with the help of the FFTW3 library \cite{FFTW05}, observing that the action of the tensor product of a discrete Fourier matrix and $I_2$ corresponds to the computation of two FFT tranforms of length $n$ on strided subvectors.
In our numerical tests, $\mathcal{K}_{\hat{S}}$ is a GMRES solver with a relative stopping tolerance $10^{-6}$.

As comparison solver we consider another preconditioning technique that does not require to assemble the Schur complement, namely the 
Least Squares Commutators (LSC) of \cite{LSC:01,LSC:06}.
It is based on the idea that one can approximate the inverse of the Schur complement, without considering the contribution of the block $E$, by
\[
	\bar{S}^{-1} = \frac{1}{\dt} \widetilde{\left(DG\right)^{-1}} DNG \widetilde{\left(DG\right)^{-1}}
	.
\]
Matrix $\overline{S}$ is never assembled, but the action of $\bar{S}^{-1}$ is computed with the above formula, where we have indicated with $\widetilde{\left(DG\right)^{-1}}$ the application of a solver for the matrix $\frac{1}{\dt}DG$, which we denote with $\mathcal{K}_{DG}$.
In our tests, we have chosen for $\mathcal{K}_{DG}$
a preconditioned conjugate gradient solver with relative stopping tolerance of $10^{-5}$, since, in the incompressible framework, the product $\frac{1}{\dt}DG$ is a Laplacian.
To provide a circulant preconditioner for $\mathcal{K}_{DG}$, it is enough to consider the block circulant matrix generated by $\dsym(\theta)\gsym(\theta)$ defined as in Remark \ref{rem:sym_dg}. Note that, for $\theta=0$, $\dsym(\theta)\gsym(\theta)$ is the null matrix, therefore in order to avoid singular matrices we introduce a rank-two correction and define the whole preconditioner for the product $\frac{1}{\dt}DG$ as 
\begin{equation}\label{eq:prec_dg}
\mathcal{P}_n:=C_n(\dsym\gsym)+\frac{1}{(2n)^2}{\bf 1}^T{\bf 1}\otimes  I_2
\end{equation}
again with ${\bf 1}=[1,\ldots,1]\in\mathbb{R}^{n}$.

For a complete Navier-Stokes simulation, the solver $\mathcal{K}_{\mathcal{A}}$ is applied at each iteration of the main non-linear Picard solver that computes a timestep. In all numerical tests, $\mathcal{K}_{\mathcal{A}}$ is a FGMRES solver with relative tolerance of $10^{-8}$.

\begin{table}[!b]
	\centering
	\small
	\begin{tabular}{r|c|c|c |c|c| c|c|c|c|}
		\toprule
		&\multicolumn{3}{c|}{$\mathcal{C}_n$} & \multicolumn{2}{c|}{$C_n(\ssym)$} & \multicolumn{4}{c|}{LSC with $\mathcal{P}_n$}\\
		$n$ &
		$\mathcal{K}_{\mathcal{A}}$ & $\mathcal{K}_{\hat{S}}$ & time (\SI{}{\second})
		& $\mathcal{K}_{\mathcal{A}}$ & $\mathcal{K}_{\hat{S}}$& 
		$\mathcal{K}_{\mathcal{A}}$ & $\mathcal{K}_{\bar{S}}$ & $\mathcal{K}_{DG}$ & time (\SI{}{\second})\\
		\midrule
		10   & 2 & 11 -- 12 & $\SI{2.61d-2}{}$ & 2 & 15 -- 16 & 2 & 2 -- 10 & 5 -- 6 & $\SI{2.08e-01}{}$\\
		20   & 2 & 10 -- 11 & $\SI{1.53e-1}{}$ & 2 & 20       & 2 & 4 -- 12 & 5 -- 6  & $\SI[retain-zero-exponent=true]{1.58e+00}{}$ \\
		40   & 2 &  9 -- 11 & $\SI{3.00e-1}{}$ & 2 & 24  	  & 2 & 3 -- 14 & 6 -- 7  & $\SI[retain-zero-exponent=true]{3.70e+00}{}$\\
		80   & 2 &  9 -- 10 & $\SI{5.57e-1}{}$ & 2 & 31  	  & 2 & 3 -- 14 & 5 -- 7  & $\SI[retain-zero-exponent=true]{7.85e+00}{}$\\
		160  & 2 &  8 -- 9 & $\SI[retain-zero-exponent=true]{1.45e+0}{}$ & 2 & no conv. 	  & 2 & 1 -- 14 & 4 -- 8  & $\SI{2.06e+01}{}$\\
		320  & 2 &  8 -- 9 & $\SI[retain-zero-exponent=true]{7.48e+0}{}$ & 2 & no conv. 	  & 3 & 1 -- 21 & 6 -- 8  & $\SI{2.42e+02}{}$\\
		640  & 2 &  7 -- 9 & $\SI{4.96e+1}{}$ & 2 & no conv. 	  & 4 & 7 -- 24 & 3 -- 10 & $\SI{2.65e+03}{}$\\
		1280 & 2 &  7 -- 8 & $\SI{3.82e+2}{}$ & 2 & no conv. 	  & 7 & 9 -- 28 & 3 -- 10 & $\SI{4.55e+04}{}$\\
		\bottomrule
	\end{tabular}
	\caption{Iterations of the solvers in the 2D pipe test with constant cross-section.
		$\mathcal{K}_{\hat{S}}$ refers to our approach, while $\mathcal{K}_{\bar{S}}$ and $\mathcal{K}_{DG}$ refer to the LSC approach. The times are the total CPU time spent in the main Krylov solver $\mathcal{K}_{\mathcal{A}}$ and its sub-solvers.}
	\label{tab:2Dpipe}
\end{table}

\begin{figure}[!t]
	\begin{center}
		\begin{tabular}{cc}
			\includegraphics[width=0.4\textwidth]{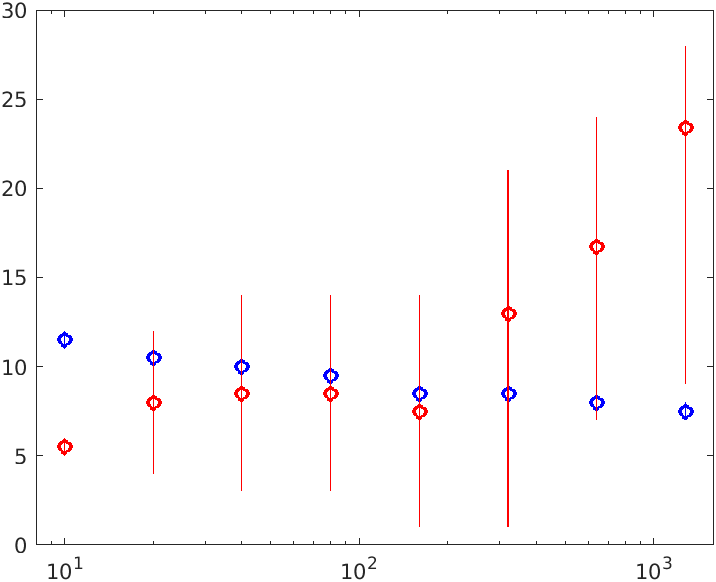}&
			\includegraphics[width=0.4\textwidth]{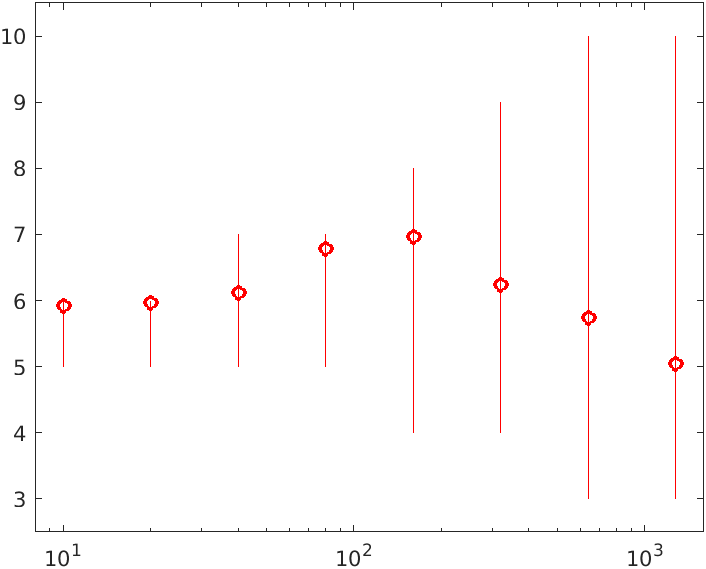}\\
			(a) & (b)
		\end{tabular}
	\end{center}
	\caption{(a) The average number and the range of iterations of $\mathcal{K}_{\hat{S}}$ in blue and of $\mathcal{K}_{\bar{S}}$ in red; (b) The average number and the range of iterations of $\mathcal{K}_{DG}$.}
	\label{fig:2Dpipeiterksp}
\end{figure}

\paragraph{Pipe with constant cross-section}

In the first test we consider a 2D pipe with constant cross-section $d(x)= \SI{0.025}{\metre}$. In inlet we impose a parabolic velocity profile with flow rate $\SI{5d-6}{\metre^2\per\second}$, while at the outlet we fix a null pressure.
Of course there would be no need to use a numerical model to compute the solution in this particular geometry, since an exact solution is known, but we conduct this as a test to verify the performance of our solver. Using $n_{x} = 1$ and $n_{y} = 3$ this setting is exactly the one adopted in \S\ref{sec:prelim} and \S\ref{sec:spectr}.
\\The main solver $\mathcal{K}_{\mathcal{A}}$ converges in at most $2$ iterations, while the number of iterations of $\mathcal{K}_{\hat{S}}$ stays constant as the number of cells grows which confirms that the block circulant preconditioner $\mathcal{C}_n$ in \eqref{eq:circ_prec} is optimal, Table~\ref{tab:2Dpipe}. For this example we also check the performances of the block circulant preconditioner $C_n(\ssym)$ in $\mathcal{K}_{\hat{S}}$. Looking again at Table~\ref{tab:2Dpipe}, we see that in this case the inner solver $\mathcal{K}_{\hat{S}}$ does not converge when the number of cells increases. The discrepancy in the performances of $C_n(\ssym)$ compared with those of $\mathcal{C}_n$ is in line with the results in Fig.~\ref{fig:symbolSwithM:fSM}(a) that clearly show how good $\ssym_{\dx}$ matches the spectrum of the Schur complement compared with $\ssym$.

Concerning the LSC approach, the number of iterations of $\mathcal{K}_{DG}$ does not grow significantly with $n$, indicating that the block circulant preconditioner $\mathcal{P}_n$ in \eqref{eq:prec_dg} for $\frac{1}{\dt}DG$ is optimal, see also Fig. \ref{fig:2Dpipeiterksp}(b). The full solver for $\mathcal{A}_n$, however, needs considerably more time to reach the required tolerance, for two reasons: 1) the number of iterations of $\mathcal{K}_{\hat{S}}$ in our approach is lower than those of  $\mathcal{K}_{\bar{S}}$ in LSC (see Fig.~\ref{fig:2Dpipeiterksp}(a)); 2) the LSC approach invokes the inner solver $\mathcal{K}_{DG}$ twice per each iteration of $\mathcal{K}_{\bar{S}}$, affecting the final computation time.

\paragraph{Pipe with variable cross-section}
In this second test we consider a 2D pipe with variable cross-section, where $d(x)$ decreases linearly from $\SI{0.025}{\meter}$ to $\SI{0.0125}{\meter}$. To perform the simulations we impose the same boundary conditions as in the previous test and again take $n_{x} = 1, n_{y} = 3$.
In Table~\ref{tab:nozzle} we compare the number of iterations computed by $\mathcal{K}_{\hat{S}}$ considering as preconditioners
\begin{enumerate}
\item $\mathcal{D}_n(\frac{1}{d}C_n(\ssym_{\dx})+\mathcal{R}_n)$, with $\mathcal{D}_n$ a diagonal matrix whose entries are an equispaced sampling of $d(x)$ on its domain (see Remark \ref{rem:general}), and $\mathcal{R}_n=\frac{1}{(2n)^2}{\bf 1}^T{\bf 1}\otimes\begin{bmatrix}1&1\\1&1\end{bmatrix}$;	
\item $\mathcal{C}_n$ with $d=\bar{d}$, that is equal to the average of the cross-section along the pipe.
\end{enumerate}

\begin{table}[!b]
	\centering
	\small
	\begin{tabular}{r|c|c|c|c|c|c|}
		\toprule
		&\multicolumn{3}{c|}{$d(x)$ in $\mathcal{K}_{\hat{S}}$} & \multicolumn{3}{c|}{$d(x) = \bar d$ in $\mathcal{K}_{\hat{S}}$}\\
	$n$ &\multicolumn{1}{p{0.12\textwidth}|}{\small non linear solver}   & $\mathcal{K}_{\mathcal{A}}$ & $\mathcal{K}_{\hat{S}}$ &  \multicolumn{1}{p{0.12\textwidth}|}{\small non linear solver}  & $\mathcal{K}_{\mathcal{A}}$ & $\mathcal{K}_{\hat{S}}$ \\
		\midrule
		10   & 6 & 1--2 & 12 -- 13 & 6 & 1--2 & 14 -- 15 \\
		20   & 5 & 1--2 & 11 -- 12 & 5 & 1--2 & 15 -- 16\\
		40   & 3 & 1--2 & 10 -- 12 & 3 & 1--2 & 14 -- 16\\
		80   & 2 & 1--2 &  9 -- 11 & 2 & 1--2 & 13 -- 16\\
		160  & 2 & 1--2 &  9 -- 11 & 2 & 1--2 & 13 -- 16\\
		320  & 2 & 1--2 &  9 -- 11 & 2 & 1--2 & 13 -- 16\\
		640  & 2 & 1--2 &  9 -- 11 & 2 & 1--2 & 13 -- 16\\
		1280 & 2 & 1--2 &  9 -- 11 & 2 & 1--2 & 13 -- 17\\
		\bottomrule
	\end{tabular}
	\caption{Iterations of the solvers in the 2D pipe test with variable cross section $d(x)$.
	In the left part, we use a diagonal scaling (defined through $d(x)$) of the block circulant preconditioner $\mathcal{C}_n$; on the right, we use $\mathcal{C}_n$ with $d=\bar{d}$, that is equal to the average of the cross-section along the pipe.}
	\label{tab:nozzle}
\end{table}

\begin{figure}[!h]
	\begin{center}
		\includegraphics[width=0.4\textwidth]{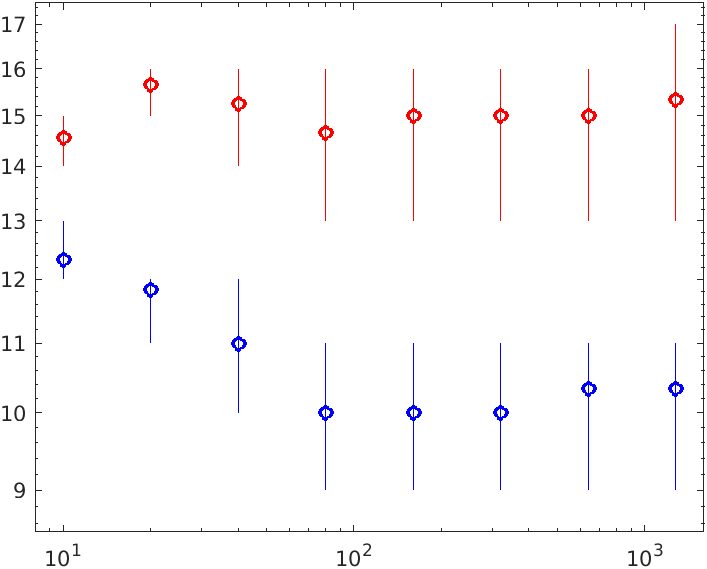}
	\end{center}
	\caption{The average number and the range of iteration of $\mathcal{K}_{\hat{S}}$ for a 2D pipe with variable cross-section. The blue values are obtained employing as preconditioner in $\mathcal{K}_{\hat{S}}$ a diagonal scaling (defined through $d(x)$) of the block circulant preconditioner $\mathcal{C}_n$; the red values are obtained using $\mathcal{C}_n$ with $d=\bar{d}$, that is equal to the average of the cross-section along the pipe.}
	\label{fig:2DnozzleiterKsp}
\end{figure}

In the first case the $\mathcal{K}_{\hat{S}}$ converges in a number of iterations that does not increase significantly with $n$, showing its optimality. Approximating the channel width with a constant value instead, avoids the diagonal matrix multiplication in the preconditioner, but causes a slightly faster increase of the iteration counts for $\mathcal{K}_{\hat{S}}$, refer to Fig.~\ref{fig:2DnozzleiterKsp}.

\paragraph{Using higher polynomial degree in the transversal direction}
In this test we analyse the efficiency of the preconditioner $\mathcal{C}_n$ in $\mathcal{K}_{\hat{S}}$ when considering different polynomial degrees $n_{y}$ in the transversal direction for the velocity, but fixed $n_x=1$ for the pressure variable. 
In this setting, we expect symbols for (1,1)-block of the coefficient matrix to take values in $\C^{2(n_y-1)\times2(n_y-1)}$, those for $(1,2)$- and $(2,1)$-blocks in
$\C^{2(n_y-1)\times2}$ and $\C^{2\times2(n_y-1)}$ respectively,
while those for the $(2,2)$-block and the Schur complement will still take values in $\C^{2\times2}$, irrespectively of $n_y$. On such basis, we can readily apply $\mathcal{C}_n$ in $\mathcal{K}_{\hat{S}}$ being sure that the sizes of all the involved matrices are consistent.

\begin{table}
	\centering
	\small
	\begin{tabular}{r|c|c|c|c|c|c|}
		\toprule
		&\multicolumn{2}{c|}{$n_{y}=4$ } & \multicolumn{2}{c|}{$n_{y}=5$}& \multicolumn{2}{c|}{$n_{y}=6$}\\
		$n$ &  $\mathcal{K}_{\mathcal{A}}$ & $\mathcal{K}_{\hat{S}}$ & $\mathcal{K}_{\mathcal{A}}$ & $\mathcal{K}_{\hat{S}}$ & $\mathcal{K}_{\mathcal{A}}$ & $\mathcal{K}_{\hat{S}}$\\
		\midrule
		10   & 2 & 10 -- 11 & 2 & 10 -- 11 & 2 & 10 -- 11 \\
		20   & 2 & 10 -- 11 & 2 & 10 -- 11 & 2 & 10 -- 11 \\
		40   & 2 & 10 -- 11 & 2 &  9 -- 11 & 2 &  9 -- 11 \\
		80   & 2 &  9 -- 10 & 2 &  9 -- 10 & 2 &  9 -- 10 \\
		160  & 2 &  8 -- 10 & 2 &  8 -- 10 & 2 &  8 -- 11 \\
		320  & 2 &  8 -- 10 & 2 &  8 -- 10 & 2 &  8 -- 11 \\
		640  & 2 &  8 -- 11 & 2 &  8 -- 11 & 2 &  8 -- 11 \\
		1280 & 2 &  8 -- 11 & 2 &  8 -- 11 & 2 &  8 -- 11 \\
		\bottomrule
	\end{tabular}
	\caption{Range of iterations for $\mathcal{K}_{\mathcal{A}}$ and $\mathcal{K}_{\hat{S}}$, in a 2D pipe with constant cross-section, with different polynomial degree in the transversal direction for the velocity.}\label{tab:2DpipeMyAlto2D}
\end{table}
Taking again the constant cross-section case, we increase $n_y$ to $4,5$ and $6$ and report the results in Table~\ref{tab:2DpipeMyAlto2D}. We note that, despite the \lq\lq looser'' approximation in the preconditioner, the solver $\mathcal{K}_{\hat{S}}$ still converges in an almost constant number of iterations when $n$ increases. From this example we can infer that the symbol of the preconditioner for the Schur complement is not changing much as far as $n_x$ stays fixed to 1.

\paragraph{3D case}
To perform a three-dimensional test, we consider a pipe 
with width equal to the 2D nozzle case above and with the same height, so that the square section area decreases quadratically
from $\SI{6.25d-4}{\square\meter}$ to $\SI{1.56d-4}{\square\meter}$. 
At the inlet we fix a constant flow rate of $\SI{5d-6}{\metre^3\per\second}$ with a parabolic profile in both the transverse directions. 

The solution is computed using different combinations of transverse polynomial degrees $n_{y}$ and $n_{z}$ for the velocity, fixed $n_x=1$ for the pressure variable.

Thanks to the matrix-sizes match pointed out in remark \ref{rem:3D}, one could be tempted to directly apply the preconditioner $\mathcal{C}_n$ in $\mathcal{K}_{\hat{S}}$ derived for the two-dimensional case also to the three-dimensional case, but results not reported here show that such choice causes high iteration numbers and sometimes stagnation of the outer nonlinear solver.
 
The reason for these poor performances may be understood by noticing that the two dimensional discretization represents in, the three dimensional setting, a flow between infinite parallel plates at a distance $d(x)$. It is not surprising that using such a flow to precondition the computation in a three dimensional pipe is not optimal. More precisely the two dimensional setting can be understood as choosing $n_{z} = 0$ in 3D. However, constant shape functions in the $z$ direction can not match the zero velocity boundary condition on the channel walls and only $n_{z} \geq 2$ would allow to satisfy them.

Fixing $n_y=3$, $n_z=2$ and following the same steps of \S\ref{sec:spectr}, we have computed an ad hoc block circulant preconditioner for the three-dimensional case. For this special choice of $n_y$ and $n_z$ the symbols of the various matrices involved in the discretization are matrix-valued with the same size as in \S\ref{sec:spectr}, but now for a fixed $n$, i.e. for a fixed $\dx$, the generating function associated with the scaled Schur complement $\frac{1}{\dt}S_n$ shows a dependency on the cross-sectional area and is given by
\begin{equation}\label{eq:3dsymbol}
	\ssym_{\dx}(\theta)  = \tfrac{Area}{c}
	\begin{bmatrix}
	-1 - (5a(\theta)-3\dx\rho)b(\theta)c             & e^{\ii\theta} - (3a(\theta)-5\dx\rho)b(\theta)c  \\
	e^{-\ii\theta} - (3a(\theta)-5\dx\rho)b(\theta)c & -1 - (5a(\theta)-3\dx\rho)b(\theta)c \\
	\end{bmatrix}	
	,
\end{equation}
where  $a(\theta) = 6\left(1-cos\theta \right) \mu c +2\dx\rho$ and $b(\theta) = \frac{175}{672}\frac{\left(1-cos\theta \right)}{a(\theta) ^2 - \dx^2 \rho^2}$.
This symbol is very similar to the one of \eqref{eq:sym_deltax}, but the different constant in the function $b(\theta)$ reflects the presence of non trivial velocity shape functions in the $z$ direction.

Therefore, we use as preconditioner in $\mathcal{K}_{\hat{S}}$ the block circulant matrix generated by $\ssym_{\dx}(\theta)$ defined as in \eqref{eq:3dsymbol} properly shifted by a rank-one block circulant matrix and scaled by a diagonal matrix whose entries are given by a sampling of the function that defines the cross-sectional area of the pipe.

\begin{table}
	\centering
	\small
	\begin{tabular}{r|c|c|c|c|c|c|c|c|c|}
		\toprule
		&\multicolumn{3}{c|}{$n_{y}=3$, $n_{z}=2$} & \multicolumn{3}{c|}{$n_{y}=3$, $n_{z}=3$}& \multicolumn{3}{c|}{$n_{y}=4$, $n_{z}=4$}\\
		$n$ & \multicolumn{1}{p{0.12\textwidth}|}{\small non linear solver}& $\mathcal{K}_{\mathcal{A}}$ & $\mathcal{K}_{\hat{S}}$ &
		\multicolumn{1}{p{0.12\textwidth}|}{\small non linear solver} & $\mathcal{K}_{\mathcal{A}}$ & $\mathcal{K}_{\hat{S}}$ &
		\multicolumn{1}{p{0.12\textwidth}|}{\small non linear solver}& $\mathcal{K}_{\mathcal{A}}$ & $\mathcal{K}_{\hat{S}}$\\
		\midrule
		10   & 13 & 1--2 & 12 -- 14 & 13 & 1--2 & 12 -- 14 & 27 & 1--2 & 11 -- 13 \\
		20   & 8  & 1--2 & 13 -- 15 & 8  & 1--2 & 13 -- 15 & 34 & 1--2 & 12 -- 14 \\
		40   & 3  & 1--2 & 13 -- 15 & 3  & 1--2 & 13 -- 15 & 37 & 1--2 & 12 -- 14 \\
		80   & 3  & 1--2 & 13 -- 16 & 3  & 1--2 & 13 -- 16 & 19 & 1--2 & 12 -- 15 \\
		160  & 2  &  2   & 13 -- 17 & 2  &  2   & 13 -- 17 &  4 & 1--2 & 12 -- 15 \\
		320  & 2  &  2   & 13 -- 17 & 2  &  2   & 13 -- 17 &  3 & 1--2 & 12 -- 15 \\
		640  & 2  &  2   & 14 -- 18 & 2  &  2   & 14 -- 18 & 2  &   2  & 12 -- 15 \\
		1280 & 2  &  2   & 14 -- 18 & 2  &  2   & 14 -- 18 & 2  & 2    & 12 -- 21 \\
		\bottomrule
	\end{tabular}
	\caption{Range of iterations for $\mathcal{K}_{\mathcal{A}}$ and $\mathcal{K}_{\hat{S}}$, in a 3D pipe with variable cross-section, with different polynomial degrees in the transversal directions for the velocity.}\label{tab:NozMyAlto3D}
\end{table}

Table \ref{tab:NozMyAlto3D} shows the range of iterations for $\mathcal{K}_{\mathcal{A}}$ and $\mathcal{K}_{\hat{S}}$. In the left part we have applied the 3D block circulant preconditioner to the corresponding simulation with $n_y=3$ and $n_z=2$. As in the two-dimensional cases, the number of iterations of $\mathcal{K}_{\hat{S}}$ does not change significantly with $n$; the nonlinear solver performs an higher number of iterations (compare with Table~\ref{tab:nozzle}) for low $n$, but they reduce fast with the increasing resolution.
In the central and right part of the table we check the performance of the 3D block circulant preconditioner corresponding to $n_y=3$ and $n_z=2$ when $n_y=n_z=3$ and $n_y=n_z=4$, respectively. As in the two-dimensional examples, for $n_y=n_z=3$, the iteration numbers stay basically unchanged, despite the fact that the preconditioner is based on $\ssym_{\dx}(\theta)$ in \eqref{eq:3dsymbol} which corresponds to a different number of degrees of freedom.
For $n_y=n_z=4$ the number of iterations of $\mathcal{K}_{\hat{S}}$ are still quite moderate, but the nonlinear solver has more problems in its convergence history. This is suggesting that the actual generating function of the Schur complement for this case departs more from the one in \eqref{eq:3dsymbol} than for the case $n_y=n_z=3$.


\section{Conclusion and perspectives}
\label{sec:conclusion}

The incompressible Navier-Stokes equations have been solved in a pipe,
using a Discontinuous Galerkin discretization over one-dimensional staggered grids. The approximation of the flow is achieved by discretization only along the pipe axis, but leveraging only on high polynomial degrees in the transverse directions.
The resulting linear systems have been studied both in terms of the associated matrix structure and in terms of the spectral features of the related coefficient matrices.  In fact, the resulting matrices are of block type, each block shows Toeplitz-like, band, and tensor structure at the same time. Using this rich matrix-theoretic information and the Toeplitz, GLT technology, a quite complete spectral analysis has been presented, with the target of designing and analyzing fast iterative solvers for the associated large linear systems. At this stage we limited ourselves to the case of block circulant preconditioners in connection with Krylov solvers: the spectral clustering at 1 has been proven and the computational counterpart has been checked in terms of constant number of iterations and in terms of the whole arithmetic cost. A rich set of numerical experiments have been presented, commented, and critically discussed. 

Of course all the facets of associated problems are very numerous and hence a lot of open problems remains. For example, the spectral analysis for more general variable coefficient 2D and 3D problems (dropping the hypothesis of elongated domain) appears achievable with the GLT theory, except for the case of variable degrees which is a real challenge. Also, more sophisticated solvers related to the Toeplitz technology, including multigrid type procedures and preconditioners can be studied for the solution of the arising saddle point problems.
All these open problems will be the subject of future investigations.

\paragraph{Acknowledgements}
All the authors are members of the INdAM research group GNCS. The work
of the first author was partly supported by the GNCS-INdAM Young
Researcher Project 2020 titled \lq\lq Numerical methods for image
restoration and cultural heritage deterioration''.

	\bibliographystyle{plain}
	\bibliography{SANSbiblio}
\end{document}